\documentclass[11pt,legno]{article}
\usepackage{booktabs} 
\usepackage{array} 
\usepackage{paralist} 
\usepackage{verbatim} 
\usepackage{amsthm}
\usepackage[table,xcdraw]{xcolor}
\usepackage[titletoc,toc,title]{appendix}
\usepackage{latexsym, amsmath, amssymb, a4, epsfig, color}
\usepackage{blindtext}
\usepackage{graphicx}
\usepackage{caption}
\usepackage{subcaption}
\usepackage[utf8]{inputenc}
\usepackage[export]{adjustbox}
\usepackage{wrapfig}
\usepackage{chngcntr}
\usepackage{apptools}
\usepackage{floatrow}
\usepackage{afterpage}
\usepackage{amssymb}

\usepackage[symbol]{footmisc}
\def\correspondingauthor{\footnote{Corresponding author.}}

\DeclareMathOperator{\arcsinh}{arcsinh}
\newcommand{\ds}{\displaystyle}

\AtAppendix{\counterwithin{definition}{subsection}}
\AtAppendix{\counterwithin{theorem}{subsection}}

\graphicspath{}

\newtheorem{theorem}{Theorem}[section]
\newtheorem{cor}[theorem]{Corollary}
\newtheorem{lemma}[theorem]{Lemma}
\newtheorem{prop}[theorem]{Proposition}
\newtheorem{definition}{Definition}

\newtheorem{remark}{Remark}

\newcommand{\pf}{\noindent {\sl Proof}. \ }
\newcommand{\p}{\partial}
\newcommand{\pd}[2]{\frac {\p #1}{\p #2}}
\newcommand{\eqnref}[1]{(\ref {#1})}

\newcommand{\Dbb}{\mathbb{D}}

\newcommand{\Nbb}{\mathbb{N}}

\newcommand{\la}{\langle}
\newcommand{\ra}{\rangle}

\newcommand{\Kcal}{\mathcal{K}}

\def\Bx{{\bf x}}
\def\By{{\bf y}}

\newcommand{\Ga}{\alpha}
\newcommand{\Gb}{\beta}

\newcommand{\Gvf}{\varphi}

\newcommand{\Gt}{\theta}

\newcommand{\Gs}{\sigma}

\newcommand{\GO}{\Omega}

\newcommand{\RR}{\mathbb{R}}
\newcommand{\CC}{\mathbb{C}}
\newcommand{\NN}{\mathbb{N}}
\newcommand{\Om}{\Omega}

\newcommand{\beq}{\begin{equation}}
\newcommand{\eeq}{\end{equation}}

\numberwithin{equation}{section}
\numberwithin{figure}{section}

\begin{document}

\title{Corner effects on the perturbation of an electric potential
\thanks{\footnotesize
This work is supported by the Korean Ministry of Science, ICT and Future Planning through NRF grant No. 2016R1A2B4014530 (to M.L.) and by the Swedish Research Council under contract 621-2014-5159 (to J.H.).}}
\author{
Doo Sung Choi\thanks{\footnotesize Department of Mathematical Sciences, Korea Advanced Institute of Science and Technology, Daejeon 34141, Korea (7john@kaist.ac.kr).} \and
Johan Helsing\thanks{\footnotesize Centre for Mathematical Sciences, Lund University, 221 00 Lund, Sweden (helsing@maths.lth.se).} \and
Mikyoung Lim\thanks{\footnotesize Department of Mathematical Sciences, Korea Advanced Institute of Science and Technology, Daejeon 34141, Korea (mklim@kaist.ac.kr).} \correspondingauthor{}}

\date{\today}
\maketitle

\begin{abstract}
  We consider the perturbation of an electric potential due to an
  insulating inclusion with corners. This perturbation is known to
  admit a multipole expansion whose coefficients are linear
  combinations of generalized polarization tensors. We define new
  geometric factors of a simple planar domain in terms of a conformal
  mapping associated with the domain. The geometric factors share properties of the generalized polarization tensors and are the
  Fourier series coefficients of a generalized external angle
  of the inclusion boundary. {Since the generalized external angle contains the Dirac delta singularity at corner points, we can determine a criteria for the
  existence of corner points on the inclusion boundary in terms of the
  geometric factors.} We illustrate and validate our results with
  numerical examples computed to a high degree of precision using integral
  equation techniques, the Nystr\"om discretization, and recursively
  compressed inverse preconditioning.
\end{abstract}



\section{Introduction}

Let $\Omega$ be a simply connected bounded domain in $\RR^2$
containing the origin and with Lipschitz boundary. We suppose that the
exterior of $\overline{\Omega}$ has unit conductivity and that
$\Omega$ is insulated. Let $h$ be a harmonic function, and consider
the following conductivity problem:
\begin{equation}
\label{condeqn}
\begin{cases}
  \ds\Delta u =0 \quad&\mbox{in } \RR^2\setminus \overline{\Omega},\\
  \ds \pd{u}{\nu}=0 &\mbox{on } \p \Omega,\\
  \ds u(\Bx)- h(\Bx) =O(|\Bx|^{-1}) &\mbox{as } |\Bx|\to \infty.
\end{cases}
\end{equation}
The background potential $h$ here is perturbed by $u$ due to the
presence of the inclusion $\Omega$. One can easily express the
perturbation $u-h$ using a boundary integral equation formulation
of~(\ref{condeqn}) involving the Neumann-Poincar\'e (NP) operator.
Furthermore, the boundary integral equation framework admits a
multipole expansion of $u-h$ whose coefficients are linear
combinations of the generalized polarization tensors (GPTs), which can
also be expressed in terms of boundary integrals.

The GPTs are a sequence of real-valued tensors that are associated with
$\Omega$ and they generalize the classical polarization tensors
\cite{polya1951isoperimetric}. They can be obtained from multistatic measurements, where a high signal-to-noise ratio is needed to acquire high-order terms \cite{ammari2014target}. 
They have been used as building blocks when solving imaging
problems for inclusions with smooth boundaries \cite{AGKLY14,AKbook}.
The GPTs contain sufficient geometric information to determine
$\Omega$ uniquely \cite{AKbook}. One can suitably approximate the
conductivity, location, and shape of an inclusion, or several
inclusions, using the first few leading terms of the GPTs by adopting
an optimization framework \cite{AGKLY14,AKLZ12}.

An inclusion with corners generally induces strong scattering close to
corner points (vertices). The understanding and application of such
corner effects is a subject of great interest. Detection methods for
inclusions, from boundary measurements, have been developed in
\cite{ikehata1999enclosing}. The gradient blow-up of the electrical
potential for a bow-tie structure, which has two closely located
domains with corners, was investigated in \cite{kang2017optimal}. It
has been shown that the spectral features of the NP operator of a
domain with corners is significantly different from that of a smooth
domain \cite{Carleman-book-16,HKL17,HMM-NJP-11,HP-ACHA-13,PP-arXiv}.
It is worth mentioning that the spectrum of the NP operator has recently drawn significant attention in relation to plasmonic resonances
\cite{ammari2013spectral,KLY,  YL17}.


This paper analyzes the effects of inclusion corners on the
perturbation of an electric potential. {We present new geometric
factors that clearly reveal the boundary information of an
inclusion, including the existence of corner points.} They satisfy mutually equivalent relations with the GPTs, so that one can compute them from the GPTs and vice versa. 
The geometric factors form an infinite sequence of complex numbers.
{The main result of this paper is that the geometric factors are actually the Fourier series coefficients of a generalized external angle function (for a precise statement, see Theorem \ref{thm:mainsigma}).}
The generalized external angle function has the Dirac delta singularity at corner points. As a consequence, we can determine when there exist corner points: the geometric
factor sequence converges to zero when the inclusion has a smooth
boundary. However, it does not converge to zero, but oscillates, if
there is any corner point on the boundary of the inclusion.
{In practice, only a finite number of leading terms of the GPTs can be obtained due to the limitation of the signal-to-noise ratio in the measurements. 
If sufficiently many terms of the GPTs are given, then the partial Fourier series sum of the generalized external angle function computed with the GPTs shows isolated high peaks at corner points, as illustrated in Section \ref{sec:experi}. 
}

Our derivation is motivated by the recent result \cite{KLL14}, where
explicit relationships between the exterior Riemann mapping coefficients
and the GPTs were derived. In addition, we consider the internal
conformal mapping, with which the geometric factors are defined and
which is obtained by reflecting the exterior Riemann mapping function
across the unit circle. We prove the Fourier relation between the
geometric factors and the generalized external angle by applying the Carath\'eodory mapping theorem on the uniform convergence of conformal mappings (see Appendix \ref{sec:cara}).

One can numerically compute the GPTs for an arbitrary domain with
corners to a high degree of precision using integral equation techniques,
the Nystr\"om discretization, and recursively compressed inverse
preconditioning (RCIP) \cite{HelsOjal08}. The geometric factors can
then be obtained from their relations with the GPTs. We present some
numerical examples to validate and visualize our results.

{The rest of the paper is organized as follows: In Section
\ref{sec:GPTsRiemann} we derive explicit connections between the GPTs
and the coefficients of the Riemann mappings and define the geometric factors. In Section
\ref{sec:curvilinear} we express the geometric factors for curvilinear polygons in terms of vertices and external angles. Section \ref{sec:genLip} is about equivalent relations between the geometric factors and the GPTs for arbitrary Lipschitz domains and the investigation of corner
effects. Section \ref{sec:experi} presents numerical examples. We
prove several recursive relations in Section \ref{sec:proof}, and we conclude with some discussion.
}

\section{Generalized polarization tensors and Riemann mappings}
\label{sec:GPTsRiemann}

\subsection{Multipole Expansion}

{We identify $\Bx=(x_1,x_2)$ in $\RR^2$ with $z=x_1+{\rm{i}}x_2\in\CC$ for
notational convenience. Let $h(\Bx)=\mbox{Re}\{H(z)\}$ with
\begin{displaymath}
H(z) = \alpha_0 +\sum_{n=1}^\infty \alpha_n z^n,\quad \alpha_n=a_n^c+{\rm i}a_n^s.
\end{displaymath}
Then, it is shown in \cite{ammari2007polarization,
  ammari2013enhancement} that the solution $u$ to \eqnref{condeqn}
satisfies $u(\Bx)=\mbox{Re}\{U(z)\}$, where $U$ is a complex analytic function in $\mathbb{C}\setminus\overline{\Om}$ such that
\beq\label{eqn:multipole}
U(z)=\alpha_0 +\sum_{n=1}^\infty \left[a_n^c\left(z^n-\sum_{m=1}^\infty \frac{\gamma_{mn}^1+\gamma_{mn}^2}{z^m}\right)+{\rm i}a_n^s\left(z^n-\sum_{m=1}^\infty\frac{\gamma_{mn}^1 -\gamma_{mn}^2}{z^m}\right)\right]
\eeq
for $|z|$ sufficiently large with}
\begin{equation}
\label{GPT12}
\begin{cases}
\ds \gamma_{kn}^1=\frac{1}{4\pi k} \left[M_{kn}^{cc}-M_{kn}^{ss}+{{\rm{i}}}(M_{kn}^{cs}+M_{kn}^{sc})\right], \\[2mm]
\ds\gamma_{kn}^2=\frac{1}{4\pi k} \left[M_{kn}^{cc}+M_{kn}^{ss}-{{\rm{i}}}(M_{kn}^{cs}-M_{kn}^{sc})\right],\quad k,n\in\mathbb{N}.
\end{cases}
\end{equation}
Here the quantities $\{M_{kn}^{\Ga\Gb}\}_{k,n\in\mathbb{N}}$ ($\Ga,
\Gb\in \{c,s\}$) are the so-called (contracted) GPTs. The zero Neumann condition on $\p \Om$ in \eqnref{condeqn} implies that \beq\label{Uconst}\Im U\,=\,\mbox{constant on }\p\Om.\eeq

The GPTs are defined in terms of boundary integrals as follows: We set polar coordinates
\begin{displaymath}
 P_n^c(\Bx) = r^n \cos n\Gt\,, \qquad
 P_n^s(\Bx) = r^n \sin n\Gt\,,
\end{displaymath}
and define
\begin{equation}
\label{gpt}
M^{\Ga\Gb}_{kn} := \int_{\p \GO} P_k^\Gb(\Bx)  (-\frac{1}{2}
I - \Kcal^*_{\p\GO})^{-1}[\nu \cdot \nabla P_n^\Ga ](\Bx) \, d\Gs(\Bx)
\end{equation}
for $k,n\in\mathbb{N}$ and $\Ga,\Gb\in \{c,s\}$. The operator
$\Kcal^*_{\p\GO}$ is the Neumann-Poincar\'e (NP) operator
\begin{equation}
\label{introkd2}
\mathcal{K}^*_{\p\GO} [\Gvf] (\Bx) =
\frac{1}{2\pi} p.v.\int_{\p\GO} \frac{\la \Bx -\By, \nu_\Bx
  \ra}{|\Bx-\By|^2} \Gvf(\By)\,d\Gs(\By)\;, \quad \Bx \in \p\GO,
\end{equation}
where $\nu_\Bx$ is the outward unit normal vector to $\p\GO$ and
$p.v.$ denotes the Cauchy principal value. It was shown in
\cite{escauriaza1992regularity,verchota1984layer} that $\lambda I-
\Kcal_{\p\Om}^*$ is invertible on $L^2_0(\p \Om)$ for $|\lambda|\geq1/2$.
See \cite{ammari2007polarization} for more properties of the NP
operator.

From \eqnref{GPT12}, one can get the $M_{kn}^{\alpha\beta}$ from the
$\gamma_{kn}^j$ and vice versa. In this sense we will refer to the
$\gamma_{kn}^j$ as GPTs as well in this paper. It is worth mentioning
that the contracted GPTs have been used in making a near-cloaking
structure \cite{ammari2013enhancement, ammari2013enhancementmaxwell} and that they can be used as
shape descriptors \cite{AGKLY14}. More applications of the
GPTs can be found in \cite{ammari2013mathematical} and references
therein.

\subsection{Two Riemann mapping functions}
\label{section:mappings}

Since $\Om$ is simply connected, thanks to the Riemann mapping theorem
there exists a unique exterior Riemann mapping $\Phi : \CC\setminus
\mathbb{D} \to \CC\setminus \Omega$ of the form
\begin{equation}
\label{Phi}
\Phi[\Om](\zeta)=C\left(
\mu_{-1}\zeta+\mu_0+\frac{\mu_1}{\zeta}+\frac{\mu_2}{\zeta^2}+\cdots\right),
\end{equation}
where $\mathbb{D}$ denotes the unit disc centered at the origin, $C>0$
is a constant, and we set $\mu_{-1}=1$. We may simply write
$\Phi(\zeta)$ when the domain is clear from the context. For $k\ge 1$,
the coefficients $\mu_k$ are invariant under translation and scaling of $\Omega$. In \cite{KLL14}, an explicit relation was
derived between the exterior Riemann mapping coefficients $\mu_k$ and
the GPTs associated with $\Om$. This means that one can compute the
$\mu_k$ from the GPTs.

In this paper, we additionally consider the internal conformal mapping
$S$ that is obtained by reflecting the exterior Riemann mapping
function across the unit circle. Remind that we assume $0\in\Om$. We will then derive formulas for the
GPTs using both the coefficients of $\Phi$ and those of $S$.

By ${\Om^r}$ we denote the reflection of $\Om$ across the unit circle, {\it i.e.},
\beq
{\Om^r}:=\Bigr\{\frac{1}{\zeta}\;|\; \zeta\in \CC\setminus\overline{\Omega}\Bigr\}\cup\{0\}.
\eeq
Note that $(\Om^r)^r=\Om$.
We also define $S[\Om^r]:\mathbb{D}\rightarrow\Om^r$ by
\beq\label{RG:conformal}
S[\Om^r](w):=
\begin{cases}
\ds\frac{1}{\Phi[\Om](\frac{1}{w})}\quad&\mbox{for }w\in \mathbb{D}\setminus\{0\},\\
\ds0\quad&\mbox{for }w=0.
\end{cases}
\eeq
We may simply write $S(w)$ when the domain is clear from the context.
Then, ${\Om^r}$ is a simply connected domain containing $0$ and $S:\mathbb{D}\rightarrow {\Om^r}$ is the interior Riemann mapping corresponding to $\Om^r$ satisfying $S(0)=0$ and $S'(0)>0$.
It is obvious from \eqnref{Phi} that $S$ admits the series expansion
\beq \label{Sseries0}
S[\Om^r](w)=\frac{1}{C}\left(b_1w+b_2w^2+\cdots\right)\quad\mbox{in } \mathbb{D}
\eeq
with some complex numbers $b_k$. Note that $b_1=1$ because $\mu_{-1}=1$.

\begin{lemma}\label{eqn:mub} The coefficients of $\Phi[\Om]$ have the
  following equivalent relation with those of $S[\Om^r]$:
\begin{equation}
\label{CauchyProduct}
\mu_{k-1} +b_{k+1}+\sum_{j=2}^k b_j \mu_{k-j}=0\,,\quad k\ge 1\,.
\end{equation}
\end{lemma}
\pf For $0<|w|<1$, we have
\begin{align*} 1&=S[\Om^r](w)\cdot
  \Phi[\Om](1/w)=\left(\sum_{j=1}^\infty b_j
  w^j\right)\left(\sum_{j=-1}^\infty \mu_kw^k\right)=1+\sum_{k=1}^\infty
  \left(\sum_{j=1}^{k+1}b_j \mu_{k-j}\right)w^k
\end{align*}
Since $b_1=\mu_{-1}=1$, this proves the lemma.\qed

\subsection{Generalized polarization tensors and Riemann mappings coefficients}
\label{sec:GPTandcoeff}

{
Let $V_1$ be the analytic function in $\mathbb{C}\setminus\overline{\Om}$ such that $\Re\{V_1(z)\}$ be the solution to \eqnref{condeqn} with $h(\mathbf{x})=\Re\{z^n\}$. Similarly, we define $V_2$ the analytic function in  $\mathbb{C}\setminus\overline{\Om}$ such that $\Re\{{\rm i}V_2(z)\}$ be the solution with $h(\mathbf{x})=\Re\{{\rm i}z^n\}$. From \eqnref{eqn:multipole} we have
\begin{align*}
\left(V_1\circ\Phi\right)(\zeta)=\Phi(\zeta)^n -
\sum_{m=1}^\infty \frac{ \gamma_{mn}^1+\gamma_{mn}^2 }{\Phi(\zeta)^m},\\
\left(V_2\circ\Phi\right)(\zeta)=\Phi(\zeta)^n -
\sum_{m=1}^\infty \frac{ \gamma_{mn}^1-\gamma_{mn}^2 }{\Phi(\zeta)^m}
\end{align*}
for sufficiently large $|\zeta|$. As discussed in \eqnref{Uconst} the zero Neumann condition in \eqnref{condeqn} implies $\Im \{V_1\circ\Phi(\zeta)\},\, \Im\{{\rm i}V_2\circ\Phi(\zeta)\}$ are constants for $|\zeta|=1$
and, thus, \beq\label{extension1}
\rm{i}\Im\{V_1\circ\Phi\}+\Im\{{\rm i}V_2\circ\Phi\}\,=\,\mbox{constant for }|\zeta|=1.
\eeq
Note that 
\begin{align}
\rm{i}\Im\{V_1\circ\Phi\}+\Im\{{\rm i}V_2\circ\Phi\}
&=\frac{1}{2}\left(V_1\circ\Phi+V_2\circ\Phi\right)(\zeta)- \frac{1}{2}\overline{\left(V_1\circ\Phi-V_2\circ\Phi\right)(\zeta)}\notag\\\label{extension2}
&=\Phi(\zeta)^n - \sum_{m=1}^\infty \frac{ \gamma_{mn}^1 }{\Phi(\zeta)^m} +\overline{\sum_{m=1}^\infty \frac{ \gamma_{mn}^2 }{\Phi(\zeta)^m}},
\end{align}
where the second equality holds for sufficiently large $|\zeta|$. This induces the following lemma, which plays an essential role in deriving relations
between the GPTs and the $\mu_k$ in \cite{KLL14}.
\begin{lemma}
\label{lemma:basic}(\cite{KLL14})
The function $V_1\circ\Phi + V_2\circ\Phi$ has an entire extension.
\end{lemma}


}

To state our result, we define two multi-index sequences
$\{\mu_{n,k}\}$ and $\{b_{n,k}\}$ ($n,k\in\NN$, $k\geq n$) such that
the following formal expansions hold:
%
\beq\notag
\sum_{k=n}^{\infty} \mu_{n,k}{x^k}=\left(\sum_{k=1}^\infty \mu_{k-2}x^k \right)^n,\quad
\sum_{k=n}^{\infty}{b_{n,k}}{x^k}=\left(\sum_{k=1}^\infty b_k x^k \right)^n.
\eeq
In other words
\begin{align*}
\mu_{n,k}
&=\sum_{s_1+s_2+\cdots+s_k =n, \atop s_1+2s_2+\cdots+ks_k=k} \frac{n!}{s_1!s_2!\cdots s_k!}  {\mu_{-1}}^{s_1}{\mu_0}^{s_2}\cdots {\mu_{k-2}}^{s_k},\\
b_{n,k}
&=\sum_{s_1+s_2+\cdots+s_k =n, \atop s_1+2s_2+\cdots+ks_k=k} \frac{n!}{s_1!s_2!\cdots s_k!} {b_1}^{s_1}{b_2}^{s_2}\cdots {b_{k}}^{s_k}.
\end{align*}
Here, $s_1,\ldots,s_k$ are non-negative integers.
In particular, we have
\begin{equation}
\label{bkk}
b_{1,k} = b_{k}\,,\quad \mu_{1,k}=\mu_{k-2}\,, \quad
\mu_{k,k}=\mu_{-1}^k=1\,,\quad b_{k,k}=b_1^k=1\,.
\end{equation}
We can deduce the following proposition using $\mu_{n,k}$ and
$b_{n,k}$. See Section \ref{proof:theorem1} for a detailed proof.
Lemma \ref{lemma:basic} plays an essential role in the derivation.

\begin{prop}\label{theorem1}
  The GPTs associated with $\Om$ have recurrence formulas with the
  coefficients of $\Phi[\Om]$ and $S[\Om^r]$.
For each $k,n\in \mathbb{N}$, we have
\begin{align}
\label{Gamma1}
\ds\gamma_{kn}^1&= C^{k+n}\left(
\mu_{n,2n+k}-\sum_{m=1}^{k-1}\frac{\gamma_{mn}^1}{C^{m+n}} b_{m,k}\right),\\[1.5mm]
\ds\gamma_{kn}^2 &=
\begin{cases}
\ds- C^{k+n}\left(\bar{\mu}_{n,2n-k}
+\sum_{m=1}^{k-1}\frac{\gamma_{mn}^2}{C^{m +n}} b_{m,k}\right),
\qquad &k\le n\,,\\[1.5mm]
\ds-C^{k+n}\sum_{m=1}^{k-1}
\frac{\gamma_{mn}^2 }{C^{m+n}} b_{m,k}\,, \qquad &k\ge n+1\,.
\label{Gamma2}
\end{cases}
\end{align}
\end{prop}
We can also express $b_k$ and
$\mu_k$ by the GPTs as follows. See Section
\ref{sec:proof:others} for a proof.
\begin{prop}\label{GPTtob}
We have
\begin{align}
\label{eqn:C}
C &= \sqrt{-\gamma_{11}^2}\,,\\
b_k& = \sum_{m=2}^{k} \frac{\gamma_{m1}^2}{C^{m+1}} b_{m,k}\,,
\quad k\geq 2\,.
\label{eqn:bk1}
\end{align}
\end{prop}
For example, with $k=2,3$ in \eqnref{eqn:bk1} we deduce, using
$b_{2,3}=2b_2$, that
$$
b_2=\gamma_{21}^2\left(-\gamma_{11}^2\right)^{-\frac{3}{2}},\quad
b_3 =2(\gamma_{21}^2)^2\left(-\gamma_{11}^2\right)^{-3}
       +{\gamma_{31}^2}\left(-\gamma_{11}^2\right)^{-2}\,.
$$
Now, with $k=1,2$ in \eqnref{CauchyProduct} we have
$
\mu_0=-b_2$ and $
\mu_1 = -b_3-b_2\mu_0$.
\begin{remark}
\label{remark:bk2}
We can show, as proved in Section \ref{sec:proof:others}, that $b_k$
also satisfies
\begin{equation}
\gamma_{11}^1 b_k = {C^{2}} \left(\mu_{k} - \sum_{m=2}^{k}
  \frac{\gamma_{m1}^1}{C^{m+1}} b_{m,k}\right), \quad k\geq 2\,.
\label{eqn:bk2}
\end{equation}
Through this relation as well as \eqnref{eqn:bk1}, we can see that
there are domain-independent relationships among the GPTs.
Additional examples are provided in \cite{KLL14}.
\end{remark}

\subsection{Geometric factor $\sigma_k$ and equivalent relations}

\begin{definition}
  We define a new sequence of geometric factors
  $\{\sigma_k\}_{k=1}^\infty$ for $\Omega$ as
\begin{equation}
\label{sigma}
  \sigma_k[\Om] := k(k+1)b_{k+1} - \sum_{j=1}^{k-1} (j+1)b_{j+1}
  {\sigma_{k-j}}\,,\qquad k\geq 1\,.
\end{equation}
\end{definition}
We can simplify this definition as
$
\ds\{\sigma_k\}_{k=1}^\infty=\mathcal{P}\big(\{b_k\}_{k=2}^\infty\big)
$
or
\beq
\label{eqn:P}
\ds \sigma_k=\sum_{i_1+2i_2+\cdots+ki_k=k}P_{k, i_1,\ldots,i_k}
b_{2}^{i_1}\cdots b_{k+1}^{i_k}\,,
\eeq
with some integer coefficients $P_{k,i_1,\ldots,i_k}$. Here,
$i_1,\ldots,i_k$ are non-negative integers.

The definition \eqnref{sigma} implies that $\mathcal{P}$ is
invertible. Thanks to the relations between the GPTs and the Riemann
mapping coefficients in the previous section, we can see that there are
mutually equivalent connections among the GPTs $\gamma_{kn}^j$, the exterior
Riemann mapping coefficients $\mu_k$, the interior Riemann mapping coefficients $b_k$, and the geometric
factors $\sigma_k$. For instance, we have
\begin{alignat*}{3}
&b_2=\frac{1}{2}\sigma_1,
\quad
&& b_3=\frac{1}{6}\left(\sigma_2 + \sigma_1^2\right),\quad
&& b_4=\frac{1}{24}\left(2\sigma_3+3\sigma_1\sigma_2+\sigma_1^3\right)\\
&\mu_0 = -\frac{1}{2}\sigma_1\,,
\quad
&&\mu_1=\frac{1}{12}\left(- 2\sigma_2 +\sigma_1^2\right)\,,
\quad&&\mu_2
 = \frac{1}{24}\left(-2\sigma_3+\sigma_1\sigma_2\right)\,,
\end{alignat*}
and
\begin{alignat*}{3}
&\gamma^1_{11}= - \frac{C^2}{6}\sigma_2 + \frac{C^2}{12}\sigma_1^2,\quad
\gamma_{12}^1=-\frac{C^3}{6}\sigma_3+\frac{C^3}{4}\sigma_1\sigma_2- \frac{C^3}{12}\sigma_1^3,\quad
 \gamma_{21}^1=-\frac{C^3}{12}\sigma_3+\frac{C^3}{8}\sigma_1\sigma_2- \frac{C^3}{24}\sigma_1^3\,,\\
 &\gamma_{11}^2=-C\bar{C}\,,\qquad
\gamma_{12}^2= C\bar{C}^2\sigma_1\,,\qquad
\gamma_{21}^2=\frac{C^2\bar{C}}{2} \sigma_1\,.
\end{alignat*}

\section{Geometric factor for a curvilinear polygon}\label{sec:curvilinear}

In this section we restrict $\Om$ to be a curvilinear polygon, in
other words $\Om=P^r$ with a simple polygon $P$, and deduce explicit connections
between the geometric factors and the corner geometry of $\Om$. Figure
\ref{figure:triangle}(a,b) illustrate a curvilinear polygon and its
reflection across the unit circle.
More specifically, we let $P$ be a simply connected region bounded by
a polygon whose vertices are $A_1,\ldots,A_n$ (ordered consecutively)
with $n\geq 3$ and with external angles $\Gb_1\pi,\ldots,\Gb_n\pi$. We
assume $0\in P$. Note that
 $-1<\beta_j<1$ and $\sum_{j=1}^n \beta_j =2$. 
 
It is well known that the interior conformal mapping of a simple
polygon can be expressed as the {\it Schwarz-Christoffel integral}.
For the polygon $P$ described above, we have
\beq\notag S[P](z) = C_1
\int_{0}^{z} {\prod_{j=1}^n \left(w - a_j \right)^{-\beta_j}}dw +
C_2\,, \eeq 
where $C_1$ and $C_2$ are complex constants, and
$a_1,\ldots,a_n$ are $n$ distinct pre-vertices on $\partial
\mathbb{D}$ satisfying $S(a_j)=A_j$ for each $j=1,\ldots,n$. One can
find a detailed explanation of the Schwarz-Christoffel integral in
many textbooks, for example in \cite{SS}.

Assuming $S(0) = 0$ and $S'(0)=1/C>0$, we have a slightly different
formulation:
\begin{equation}
S[P](z)=\frac{1}{C}\int_{0}^{z} \pi(w)\,dw\,,\qquad
\pi(z)={\prod_{j=1}^n \left(1 - \frac{z}{a_j}\right)^{-\beta_j}}\,.
\end{equation}
Here, the branch-cut is given such that $\pi(0)=1$. The coefficients
in the expansion \eqnref{Sseries0} satisfy $C\cdot
S'(z)=1+\sum_{k=1}^\infty b_{k+1} (k+1) z^{k}$. From the fact
$\pi(z)=C\cdot S'(z)$ we deduce \beq\label{piz}
b_{k+1}=\frac{\pi^{(k)}(0)}{(k+1)!}\,.\eeq

\begin{lemma}
\label{lemma:sigma}
For the polygon $P$ described above, we have for each $k\in\NN$
\beq\label{def:sigma}\sigma_k[P^r]=\sum_{j=1}^n \beta_j a_j^{-k}.\eeq
\end{lemma}
\pf
Set $\widetilde{\sigma}_k=\sum_{j=1}^n \beta_j a_j^{-k}$ and consider the function
\beq\notag
F_k(w)=(k-1)!\sum_{j=1}^n \beta_j{a_j}^{-k} \left(1-\frac{w}{a_j} \right)^{-k},\quad k\geq 1\,.\eeq
One can easily see that $\pi'=\pi F_1$ and ${F_k}'= F_{k+1}$. Applying the Leibniz rule, we have
\begin{equation*}
\pi^{(k)} = (\pi F_1)^{(k-1)} = \sum _{j=0}^{k-1}{k-1 \choose j}{\pi}^{(j)}{F_1}^{(k-1-j)}
 = \sum _{j=0}^{k-1}{k-1 \choose j}{\pi}^{(j)}{F_{k-j}}.
\end{equation*}
Note that $F_k(0)=(k-1)!\,\widetilde{\sigma}_k$. Evaluating the above equation at $0$, we prove
\begin{equation}
\pi^{(k)}(0) = (k-1)!\,\widetilde{\sigma}_k + \sum_{j=1}^{k-1} \frac{(k-1)!}{j!}\pi^{(j)}(0){\widetilde{\sigma}_{k-j}}.\label{eqn:pisigma}
\end{equation}
The lemma then follows as a direct consequence of \eqnref{piz} and \eqnref{sigma}.\qed

\section{Analysis of corner effects}\label{sec:genLip}

We now consider an arbitrary planar domain with corners.  We assume that
$\Omega$ is a simply connected domain bounded by a piecewise regular analytic curve with a finite number of corners. Figures
\ref{fig:asymm}(a,b) illustrate an example of
such domain and its reflection across the unit circle.

We characterize the corner effects in the geometric factor (and as a
result, in the GPTs or in Riemann mapping coefficients) by
approximating $\Om^r$ with simple polygons.

\subsection{Generalized external angle of $\p {\Om^r}$ }

Let $S[\Om^r]:\mathbb{D}\rightarrow \Om^r$ be given as in Section \ref{section:mappings}. Since
$\partial \Omega$ is a piecewise regular analytic curve, so is $\partial
{\Om^r}$. Since $\partial {\Om^r}$ is a Jordan curve, the interior Riemann
mapping $S$ extends to a bijective continuous function from
$\overline{\mathbb{D}}$ to $\overline{{\Om^r}}$. Let $\zeta$ be an
arbitrary point on $\p \Om^r$ which is not a corner point. Then the
inverse mapping $S^{-1}:\Om^r \rightarrow\mathbb{D}$ extends
holomorphically across $\zeta$ (see \cite{cho1996regularity}) and
$(S^{-1})'(\zeta)\neq 0$ (see \cite[Theorem 18, p. 217, and Theorem
20, p. 226]{clebsch1931mathematische}). From the inverse function
theorem, $S$ is smooth near $S^{-1}(\zeta)$. We assume that the curvature of $\p
\Om^r$ is uniformly bounded except at the corner points.

Denote $\alpha(t)=(x(t),y(t))$, $t\in[0,1]$, the piecewise analytic
parametrization of $\p \Om^r$ by the Riemann mapping $S[\Om^r]$. In other
words
\begin{displaymath}
\alpha(t)=S[\Om^r](e^{2\pi t{\rm i}})\,,\qquad t\in[0,1]\,,
\end{displaymath}
together with finitely many points $0\leq t_1<\cdots<t_M<1$, where $\alpha(t_l)$, $l=1,\ldots,M$, are corner points on $\p\Om^r$. We set
$t_{M+1}=t_1+1$ and regard $\alpha$ as a periodic function with period
$1$, for notational convenience. Owing to the assumption $C>0$ in
series expansions for $S$, $\alpha$ has positive orientation.  We
define the external angle $\beta_l \pi$ at each corner $\alpha(t_l)$
by the signed angle between the two vectors $\alpha'(t_l-)$ and
$\alpha'(t_l+)$, {\it i.e.}, $\beta_l\in(-1,1)$ and
\beq\label{def:beta}\beta_l \pi =
\arg\left(\alpha'(t_l+)\right)-\arg\left(\alpha'(t_l-)\right)\eeq with
a modulus of $2\pi$. We then generalize the concept of external angle
for any boundary points as follows. The generalized external angle is
essential in understanding corner effects in perturbations of an electric potential.

\begin{definition}
  Define a function $\Theta : [0,1] \to \mathbb{R}$ as
  \beq\label{def:theta}\Theta(t): = \frac{1}{\pi} k_{\rm{g}}(t) |\alpha'(t)|+
  \sum_{l=1}^M \beta_l \delta(t-t_l)\,, \quad t\in[0,1]\,,\eeq where
  $k_{\rm{g}}$ denotes the (geodesic) curvature of $\alpha$ in
  $\mathbb{R}^2$, i.e.,
  $$k_{\rm{g}}(t)=\frac{x'(t)y''(t)-x''(t)y'(t)}{|\alpha'(t)|^3}\,.$$
  We call
  $\pi\Theta(t)$ the generalized external angle of $\p {\Om^r}$ at
  $\alpha(t)$. We may indicate the associated domain by writing
  $\Theta[\p {\Om^r}]$ if necessary.
\end{definition}

Recall that the internal angle is $\pi(1-\beta)$ at the corner
$\alpha(t_l)$. From \cite{warschawski1955theorem} we have
$S^{-1}(z)\sim (z-\alpha(t_l))^\frac{1}{1-\beta}$ for $z\in\Om^r$ near
$\alpha(t_l)$ and, thus, $\alpha(t)\sim (t-t_l)^{1-\beta}$. Here, $f\sim g$ near $z_0$ means $\lim_{z\rightarrow z_0}\frac{f(z)}{g(z)}\in\mathbb{C}\setminus\{0\}$.  Therefore,
$
\frac{1}{\pi} k_{\rm{g}}(t) |\alpha'(t)|$ is integrable. From the Gauss-Bonnet formula we have \beq \int_{0}^{1} \Big[k_{\rm{g}}(t)
|\alpha'(t) | + \sum_{l=1}^N \beta_l \pi \delta(t-t_l)\Big] dt =
\int_{\partial \Omega} k_{\rm{g}}(s)ds + \sum_{l=1}^N \beta_l \pi =2\pi\,.
\label{gauss} \eeq

\subsection{Approximation of ${\Om^r}$ by a sequence of polygons}

Fix $n\in\NN$ and consider equally distanced nodes
$\big\{\frac{j}{n}|1\leq j\leq n\big\}$ on $[0,1]$. For each
$l=1,\ldots,M$, let $j_l$ be the index such that
$\alpha(\frac{j_l}{n})$ is the nearest point in
$\big\{\alpha(\frac{j}{n})|1\leq j\leq n\big\}$ to the corner point
$\alpha(t_l)$.  We set
\begin{align*}
  p_{n,j} &=\ds\begin{cases}
    \alpha(t_l),\quad &\mbox{if }j=j_l\mbox{ for some }l=1,\ldots,M,\\
    \alpha(\frac{j}{n}),\quad&\mbox{otherwise},
\end{cases}\\
e_{n,j} &=\ds
\begin{cases}
e^{2\pi t_l{\rm i}},\quad &\mbox{if }j=j_l\mbox{ for some }l=1,\ldots,M,\\
e^{2\pi\frac{ j}{n}{\rm i}},\quad&\mbox{otherwise}.
\end{cases}
\end{align*}

We denote by $P_n$ the polygon defined by the chain of edges
$[p_{n,j}, p_{n,j+1}]$, $j=1,\ldots,n$, with $p_{n,n+1}=p_{n,1}$. We
let $\beta_{n,j}\pi$ be the external angle of $P_n$ at each $p_{n,j}$.
Recall that $\beta_{n,j}\in(-1,1)$.
\begin{figure}[ht]
  \centering \includegraphics[scale=0.35]{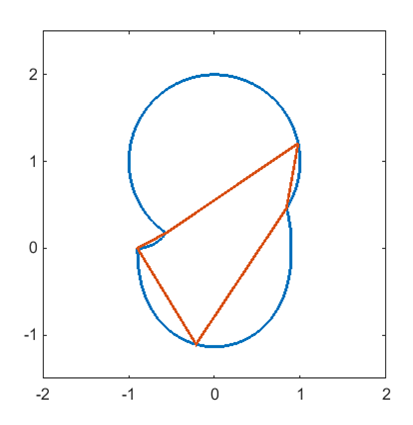}\hskip
  3mm\includegraphics[scale=0.35]{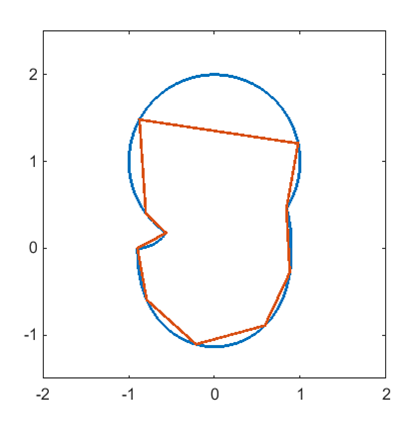}\hskip
  3mm\includegraphics[scale=0.35]{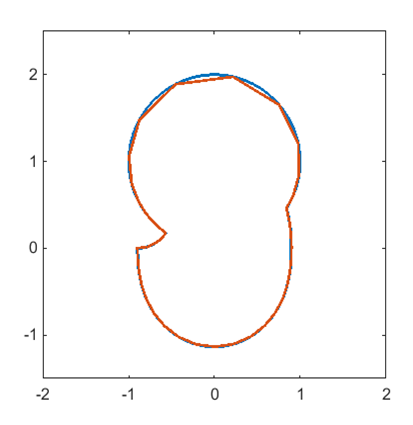}
\caption{The polygon $P_n$ (in red) approximates ${\Om^r}$ (in blue), where $\Om$ is given as in Figure \ref{fig:asymm}. The number of edges are  $n=5,10,40$ from the left to the right.}
\label{approx}
\end{figure}
The polygon $P_n$ is an $n$-sided simple polygon with sufficiently
large $n$, and one can show that the sequence of polygons $P_n$
converges in the sense of Carath\'eodory to its kernel $\Om^r$.  See
Appendix A for the definition of such convergence. Figure \ref{approx}
shows how $P_n$ approximates to $\Om^r$ as $n$ increases.  We can
obtain the asymptotic for the external angles.
\begin{lemma}
\label{lemma:betadecay}
Let $n$ be large enough so that $P_n$ is a simple polygon. For $p_{n,j}$ located away from corner points, we have \beq\label{betanj} \beta_{n,j}=
\ds\frac{1}{n\pi}{k_{\rm{g}}\left(\frac{j}{n}\right)}\left|\alpha'\left(\frac{j}{n}\right)\right|+o(\frac{1}{n}).
\eeq
\end{lemma}
\pf
We set $h=\frac{1}{n}$ and $t=\frac{j}{n}$. Note that the external
angle $\beta_{n,j}\pi$ is the relative angle of
$\alpha(t+h)-\alpha(t)$ with respect to the direction of
$\overrightarrow{\alpha(t-h)\alpha(t)}$. Let us denote
$$\Delta^{\pm}_{h} \alpha (t) = \pm \left(\alpha(t\pm h)-\alpha(t)\right).$$
Then the angle $\beta_{n,j}$ satisfies
$$\left|\beta_{n,j}\right|= \frac{1}{\pi} \arccos \left(\frac{\langle \Delta^+_{h}\alpha(t),\ \Delta^-_{h}\alpha(t) \rangle}{|\Delta^+_{h}\alpha(t) | \cdot|\Delta^-_{h}\alpha(t) |}\right).$$

Let us suppose that $p_{n,j}=\alpha(\frac{j}{n})$ is located away from corner points.
 By applying the Taylor's Theorem, we have
$$
\Delta^\pm_{h}\alpha(t) = \alpha'(t) h\pm\frac{\alpha''(t)}{2!}h^2+\frac{\alpha'''(t)}{3!}h^3+o(h^3)
$$
and, therefore,
\begin{align*}
\frac{\langle \Delta^+_{h}\alpha(t),\ \Delta^-_{h}\alpha(t) \rangle}{|\Delta^+_{h}\alpha(t) | \cdot| \Delta^-_{h}\alpha(t) |}&=1-\frac{h^2}{2}\left( \frac{x'(t)y''(t) -  x''(t)y'(t)}{[x'(t)]^2 + [y'(t)]^2} \right)^2+o(h^2)\\
&= 1-\frac{h^2}{2}k_{\rm{g}}(t)^2 | \alpha'(t) |^2 +o(h^2).
 \end{align*}
Therefore, we get
$$
\left|\beta_{n,j}\right|
 = \frac{h}{\pi} |k_{\rm{g}}(t)||\alpha'(t)|+o(h) = \ds\frac{1}{n\pi}{\left|k_{\rm{g}}\left(\frac{j}{n}\right)\right|}\left|\alpha'\left(\frac{j}{n}\right)\right|+o\left(\frac{1}{n}\right).
$$
Since the magnitude of $\beta_{n,j}$ is small due to $t$ not being a
corner point, the sign of $\beta_{n,j}$ coincides with that of the cross product 
$(\nabla_h^-\alpha(t))\times(\nabla_h^+\alpha(t))$. Since $(\nabla_h^-\alpha(t))\times(\nabla_h^+\alpha(t))$ has the same sign as $k_{\rm{g}}(\frac{j}{n})$, we prove the lemma. \qed

Since $\p\Om^r$ is a piecewise analytic curve, we have from a similar analysis as in the previous lemma that for any $\epsilon>0$, there exists $\delta=\delta(\epsilon)>0$ such that 
\beq\label{nearcorner1}
\sum_{p_{n,j}\in B_\delta}|\beta_{n,j}|<\epsilon
\eeq
with $B_\delta=\cup_{l=1}^M\left \{z\in\p P_n:0<|z-\alpha(t_l)|<\delta\right\}$.
Since $\alpha(t)\sim(t-t_l)^{1-\beta}$, we have 
\beq\label{nearcorner2}
\left|\frac{j}{n}-t_l\right|=O(\delta^{\frac{1}{1-\beta}})\quad\mbox{for }p_{n,j}\in B_\delta.
\eeq

\noindent{\textbf{Riemann mapping functions.}} 
We set $\Omega_n =(P_n)^r$ and, thus, $(\Om_n)^r = P_n$.
The Riemann mapping functions $\Phi[\Om_n]$ and $S[P_n]$ admit the series expansions
\begin{align*}
   \Phi[\Om_n](z)&=C_n\left(\mu_{n,-1}z+\mu_{n,0}+\frac{\mu_{n,1}}{z}
        +\frac{\mu_{n,2}}{z^2}+\cdots\right),
\quad z\in\mathbb{C}\setminus \mathbb{D}\,,\\
S[P_n](z)&=\frac{1}{C_n}\left(b_{n,1}z+b_{n,2} z^2 +
b_{n,3}z^3+\cdots\right),\quad z\in\mathbb{D}\,,
\end{align*}
with $\mu_{n,-1}=b_{n,1}=1$ and some constants $C_n>0$.

Since $P_n$ converges in the sense of Carath\'eodory to its kernel
${\Om^r}$, we have uniform convergence from $S[P_n]$ to $S[\Om^r]$ and from
$\Phi[\Om_n]$ to $\Phi[\Om]$ as $n\rightarrow\infty$ from the Carath\'eodory's
mapping theorem (see Appendix A for the statement and references).
Therefore, for each $k$ we have
\begin{displaymath}
\mu_{n,k}\rightarrow\mu_k\,,\quad
b_{n,k}\rightarrow b_k\,,\quad
C_n\rightarrow C\,,\quad\mbox{as }n\rightarrow\infty\,.
\end{displaymath}
Here, $\mu_k,\,b_k,\, C$ are coefficients corresponding to $\Om$.
From \eqnref{eqn:P} we have for each $k$
\beq\label{eqn:sigmaconv}
\sigma_{k}[\Om_n]\rightarrow\sigma_k[\Om]\,,\quad\mbox{as }n\rightarrow\infty\,.
\eeq

Since the boundary of $\Om^r$ is a Jordan curve, the corresponding
Riemann mapping $S[\Om^r]$ extends to a bijective continuous function from
$\overline{\mathbb{D}}$ onto $\overline{\Om^r}$ as explained before. From the Carath\'eodory's mapping theorem, $S[P_n]^{-1}$ converges uniformly to
$S[\Om]^{-1}$ on any compact subset of $\Om$ as $n \to \infty$.  Moreover,
it was proved in \cite[pp. 75--79, vol. 2]{markushevich1977theory}
that $S[P_n]^{-1}$ converges uniformly to $S[\Om]^{-1}$ on $\overline{\Om^r}$
if $P_n$ decreases to $\Om^r$. In the proof, the equicontinuity of
$\{S[P_n]^{-1}\}$ plays an essential role. By slightly modifying the
proof of \cite[Theorem 2.26, vol.3]{markushevich1977theory}, we have the following:
\begin{lemma}\label{lemma:equi}
$\left\{S[P_n]^{-1}|_{\overline{P_n\cap \Om^r}}\right\}_{n\in \NN}$
is equicontinuous.
\end{lemma}
\pf
We set $f_n = S[P_n]^{-1}$ and $f=S[\Om^r]^{-1}$. 
If the lemma is not true, then there exists $\epsilon_0>0$, a sequence $\{n_k\}$ with $n_1<n_2<\dots$, and two sequences $\{ z_k' \}$, $\{z_k''\}$ such that $z_k',z_k''\in P_n\cap \Om^r$ and
$$|f_{n_k}(z_k')-f_{n_k}(z_k'')|\geq \epsilon_0\mbox{ for each }k,\quad\mbox{while }\lim_{k\rightarrow \infty}(z_k'-z_k'')=0.$$
Due to the uniform convergence of $f_n$ to $f$ on a compact subset of $\Om^r$, we have (by taking a subsequence of $\{n_k\}$ if necessary)
\beq\notag
\lim_{k\rightarrow\infty} z_k'=\lim_{k\rightarrow\infty}z_k''=\xi\in\p
\Om^r
\eeq
and
\begin{equation}
\label{eqn:fnk}
f_{n_k}(z_k')\rightarrow w'\,,\quad
f_{n_k}(z_k'')\rightarrow w''\,,\quad
|w'|=|w''|=1\,, \quad
w'\neq w''\,.
\end{equation}

We can take a sequence $\rho_k$ such that $z_k',z_k''\in
\left\{z:|z-\xi|<\rho_k\right\}$ and $\rho_k\rightarrow 0$ as
$k\rightarrow\infty$.  From the construction of $P_n$, there exists
$R>0$ such that $\left\{z:|z-\xi|=\rho \right\}\cap P_{n_k}$ is an arc with the
center $\xi$ for any $\rho\in (\rho_k,R)$.
From \eqnref{eqn:fnk} we can assume
\begin{equation}
\label{eqn:fnk2}
|f_{n_k}(z_k')-f_{n_k}(z_k'')|>0\,,\quad
|f_{n_k}(z_k')|>0.5\,, \quad
|f_{n_k}(z_k'')|>0.5\,.
\end{equation}
Let $l_k'$, $l_k''$ be the line segments joining $0$ to
$f_{n_k}(z_k')$, $f_{n_k}(z_k'')$, and consider the two curves
$f_{n_k}^{-1}(l_k')$ and $f_{n_k}^{-1}(l_k'')$. For each
$\rho\in(\rho_k, R)$, there is an arc $\Lambda_{k,\rho}$ contained in $P_{n_k}$ with the
center $\xi$ such that one boundary point, say $z_{k,\rho}'$, is in
$f_{n_k}^{-1}(l_k')$ and the other, say $z_{k,\rho}''$, is in
$f_{n_k}^{-1}(l_k'')$.  
\smallskip

From $f(0)=0$ and the continuity of $f^{-1}$, for given $\epsilon>0$ we have $\{|z|<\delta\}\subset f\left(\{|z|<\epsilon\}\right)$ with some $\delta=\delta(\epsilon)>0$. 
Due to $|\xi|=1$, we can take $R$ sufficiently small such that
$f(\Lambda_{k,\rho})$ is located away from $0$.
From the uniform convergence of $f_n^{-1}$ to $f^{-1}$ near $0$, there is a $d>0$ independent of $k$ such
that
\beq\label{eqn:inf}\inf_{\rho\in(\rho_k,R)}\Bigr\{|f_{n_k}(z_{k,\rho}')|,
|f_{n_k}(z_{k,\rho}'')|\Bigr\}> d\,.\eeq
Therefore, from the fact $f_{n_k}(z_{k,\rho}')\in l_k'$ and $f_{n_k}(z_{k,\rho}'')\in l_k''$, there is a $\tilde{d}>0$
independent of $k$ such that
$$0<\tilde{d}<\left|f_{n_k}(z_{k,\rho}')-f_{n_k}(z_{k,\rho}'')\right|\quad\mbox{for all }\rho\in(\rho_k, R).$$
We now compute
$$0<\tilde{d}<\left|\int^{z_{k,\rho}''}_{z_{k,\rho}'}f'_{n_k}(z)dz\right|\leq\int_{\Lambda_{k,\rho}}\left|f'_{n_k}(\xi+\rho
  e^{{\rm i}\theta})\right|\rho d\theta.$$
By using the Cauchy-Schwarz inequality,
$${\tilde{d}}^2\leq 2\pi \int_{\Lambda_{k,\rho}}\left|f'_{n_k}
  (\xi+\rho e^{{\rm i}\theta})\right|^2\rho^2 d\theta\quad\mbox{for all }
\rho\in(\rho_k, R)\,.$$
By dividing both sides by $\rho$ and
integrating them, we finally have
\begin{align*}
  \tilde{d}^2\ln \frac{R}{\rho_k}& \leq 2\pi\int_{\rho_k}^R
  \int_{\Lambda_{k,\rho}}
  \left|f'_{n_k}(\xi+\rho e^{{\rm i}\theta})\right|^2\rho d\theta d\rho\\
  &\leq 2\pi\int_{B_R(\xi)\cap P_{n_k}} \left|f'_{n_k}(\xi+\rho
    e^{{\rm i}\theta})\right|^2\rho d\theta d\rho \leq
  2\pi\cdot\mbox{area}\left(f_{n_k}(P_{n_k})\right) \leq
  2\pi\cdot\mbox{area}(\mathbb{D}).
\end{align*}
The right-hand side is bounded independently of $k$. This fact contradicts $\rho_k\rightarrow 0$ as $k\rightarrow\infty$.\qed

\smallskip

 \begin{cor}\label{lem:SnS} we have
$$\sup_{j=1,\ldots,n} \left|  \left(S[P_n]^{-1}\circ S[\Om^r]\right)(e_{n,j}) -e_{n,j}\right|\rightarrow 0\quad\mbox{as }n\to \infty\,.$$\end{cor}
\pf
We set $f_n = S[P_n]^{-1}$ and $f=S[\Om^r]^{-1}$ as in the proof of the previous lemma.

For $z=p_{n,j}$, which is in $\p (P_n\cap \Om^r)$, we decompose
$$ \bigr| f_n(z) - f(z) \bigr|  \le \bigr| f_n(z) - f_n(\zeta) \bigr|  + \bigr| f_n(\zeta) - f(\zeta) \bigr|  + \bigr| f(\zeta) - f(z) \bigr|$$
with $\zeta\in P_n\cap \Om^r$ closely located to $z$.
From Lemma \ref{lemma:equi}, the uniform convergence of $f_n$ to $f$ on a compact subset of $\Om^r$ and the uniform continuity of $f$ on $\overline{\Om^r}$,  we can easily derive that
$\bigr|  f_n(z) -f(z)\bigr|\rightarrow 0$ uniformly for $j=1,\dots,n$ as $n\to \infty$.
Due to $p_{n,j}=S[\Om^r](e_{n,j})$, we finish the proof.
\qed

\subsection{Corner effects}

Since $P_n$ is a polygon, $S[P_n]$ admits \beq \notag S[P_n](z) =
\frac{1}{C_n} \int_{0}^{z} {\prod_{j=1}^n \left(1 -
    \frac{w}{a_{n,j}}\right)^{-\beta_{n,j}}}dw \eeq with a positive
constant $C_n$ and pre-vertices
$a_{n,j}=S[P_n]^{-1}(p_{n,j})=\left(S[P_n]^{-1}\circ S[\Om^r]\right)(e_{n,j})$. We denote the
geometric factor of $\Om_n$ by $\sigma_{n,k}$, in other words $\sigma_{n,k}=\sigma_k[\Om_n]$. Then, it follows from
Lemma \ref{lemma:sigma}
\begin{equation}
\label{sigma_nk}
\sigma_{n,k}=\sum_{j=1}^n \beta_{n,j} a_{n,j}^{-k},\quad k\in\NN\,.
\end{equation}

In order to have the value of $\sigma_{n,k}$ we need to know
$a_{n,j}=\left(S[P_n]^{-1}\circ S[\Om^r]\right)(e_{n,j})$, which are the pre-vertices of
$P_n$. However, the problem of finding the pre-vertices for a given
polygon, the \textit{Schwarz-Christoffel Parameter Problem}, is
challenging to solve for arbitrary polygons. There are numerical
algorithms for finding pre-vertices for special polygons, for instance
\cite{DT}. Because of this difficulty we instead consider
$$
\widetilde{\sigma}_{n,k}:=\sum_{j=1}^n \beta_{n,j} e_{n,j}^{-k}.$$

\begin{lemma}
\label{eqn:twosigma}
We have
$$\sigma_k[\Om]=\lim_{n\rightarrow\infty}\widetilde{\sigma}_{n,k}\,,
\quad k\in\NN\,.$$
\end{lemma}
\pf
Note that
$$
\widetilde{\sigma}_{n,k}
=\sum_{1\leq j\leq n,\atop j\neq j_1,\ldots,j_M}\beta_{n,j} e_{n,j}^{-k}+\sum_{j=j_1,\ldots,j_M}\beta_{n,j} e_{n,j}^{-k}\,
$$
and $|a^{-k}-b^{-k}|=|b^k-a^k|/|a^kb^k|\leq k|a-b|$ for $a,b\in\p\mathbb{D}$. Therefore, we have
\begin{align*}|\sigma_{n,k}-\tilde{\sigma}_{n,k}|
&\leq \sum_{j=1}^n k|\beta_{n,j}||a_{n,j}-e_{n,j}|\\
&=\sum_{1\leq j\leq n,\atop j\neq j_1,\ldots,j_M}k|\beta_{n,j}|
|a_{n,j}-e_{n,j}|+\sum_{j=j_1,\ldots,j_M}k|\beta_{n,j}||a_{n,j}-e_{n,j}|\,.
\end{align*}

Fix $\epsilon>0$ and choose $\delta>0$ such that \eqnref{nearcorner1} holds. 
Using Lemma \ref{lemma:betadecay} and the fact $\beta_{n,j}\in(-1,1)$,
we obtain
\beq\notag
\left|\sigma_{n,k}-\tilde{\sigma}_{n,k}\right|
\leq k \sup_{j=1,\ldots,n} \left|a_{n,j} -e_{n,j}\right|
\left(\sum_{ p_{n,j}\notin B_\delta \atop j\neq j_1,\ldots,j_M}\frac{1}{n\pi}
\left|k_{\rm{g}}\left(\frac{j}{n}\right)\alpha'\left(\frac{j}{n}\right)\right|
+\sum_{ p_{n,j}\in B_\delta }\beta_{n,j}+M\right).\eeq
Remind that $\frac{1}{\pi} k_{\rm{g}}(t) |\alpha'(t)|$ is integrable. From Corollary \ref{lem:SnS} and \eqnref{nearcorner1}, it follows for each $k$ that \beq
\left|\sigma_{n,k}-\tilde{\sigma}_{n,k}\right|\rightarrow 0\quad\mbox{as } n\rightarrow \infty\,.\notag
\eeq
This proves the lemma thanks to \eqnref{eqn:sigmaconv}.\qed

The following proposition shows that the limit is indeed the Fourier transform of $\Theta$.
\begin{theorem}\label{thm:mainsigma}
  We assume that $\Omega$ is a simply connected domain that is bounded
  by a piecewise regular analytic curve with a finite number of
  corners. Then, we have
\begin{equation}
\label{thmeqn1}
  \sigma_{k}[\Om]=\widehat{\Theta}(k)\,,\quad k\in\NN\,,
\end{equation}
  where $\widehat{\Theta}$ is the Fourier coefficients of
  $\Theta[\p\Om^r]$, that is $\widehat{\Theta}(k)=\int_{0}^{1}
  \Theta[\p\Om^r] (t) e^{-2\pi{\rm i}kt} dt$.
 \end{theorem}
\pf
{Let $\epsilon>0$ and choose $\delta>0$ such that \eqnref{nearcorner1} holds. Using Lemma \ref{lemma:betadecay} and \eqnref{nearcorner2}, we have
\begin{align*}
  \lim_{n \to \infty} \widetilde{\sigma}_{n,k}
  & = \lim_{n \to \infty} \sum_{j=1}^n \beta_{n,j}e_{n,j}^{-k} \nonumber\\
  &=\lim_{n\rightarrow\infty}\sum_{p_{n,j}\notin B_\delta \atop j\neq j_1,\ldots,j_M} \frac{1}{n\pi}k_{\rm{g}}\left(\frac{j}{n}\right)
  \left|\alpha'\left(\frac{j}{n}\right)\right|e_{n,j}^{-k}+O(\epsilon)
  +\sum_{l=1}^M \beta_l e_{n,j_l}^{-k}\\
  &=\int_{\delta^{\frac{1}{1-\beta}}}^1\frac{1}{\pi}k_{\rm{g}}(t)|\alpha'(t)|e^{-2\pi{\rm i}kt}dt
  +\sum_{j=1}^M \beta_l e^{-2\pi{\rm i}kt_l }+O(\epsilon).
  \end{align*}
  }
  Since $\epsilon$ can be arbitrarily small, we have $$\lim_{n \to \infty} \widetilde{\sigma}_{n,k}
   = \int_{0}^{1} \Theta \left(t\right) e^{-2\pi{\rm i}kt} dt
  =\widehat{\Theta}(k).$$
Thus, we prove the theorem by using Lemma \ref{eqn:twosigma}.  \qed

\smallskip
\smallskip

Since $\frac{1}{\pi} k_{\rm{g}}(t) |\alpha'(t)|$ is integrable, its Fourier
coefficient decays to zero by the Riemann-Lebesgue Lemma.  As a direct
consequence of Theorem \ref{thm:mainsigma} and the definition of
$\Theta$ we have the following criterion for the existence of corner
points:
\begin{cor}
\label{lemma:sigmadecay}
If $\partial \Omega$ has at least one corner point, then
$\sigma_k[\Om] = O(1)$ with no decay, or $\sigma_k[\Om]$ oscillates between
some bounded values in $\mathbb{C}$ as $k$ goes to infinity. Otherwise, if $\partial \Omega$
is a regular analytic curve, then $\sigma_k = O(r^k)$ with some constant $r\in (0,1)$.
\end{cor}

From \eqnref{gauss} the constant coefficient of $\Theta$ is $2$.
Hence, the Fourier series of $\Theta : [0,1] \to \mathbb{R}$ is
$$
\Theta[\p\Om^r](t) = 2 + 2\sum_{k=1}^\infty \Re\{\sigma_k\} \cos(2 \pi kt) + 2\sum_{k=1}^\infty \Im\{\sigma_k\} \sin(2 \pi kt)\,.
$$

\begin{cor}\label{radialsymm}
  $\Om$ is $n$-point radially symmetric if and only if $\sigma_k = 0$
  for every $k \not\equiv 0 \pmod{n}$.
\end{cor}
\pf If $\Om$ is $n$-point radially symmetric, so is $\Om^r$. That
means $\Theta[\p\Om^r](t) = \Theta[\p\Om^r](t+\frac{1}{n})$ for all
$t$. This is equivalent to $\Re\{\sigma_k\}=\Im\{\sigma_k\}=0$ for
every $k \not\equiv 0 \pmod{n}$. \qed

\subsection{Imaging from a finite number of components of the GPTs}\label{sec:finiteGPTs}
The GPTs can be obtained from multistatic measurements \cite{ammari2014target}.
We need infinitely many components of the GPTs to have the full sequence of geometric factors. However, one can accurately acquire only a finite number of components of the GPTs from far-field measurements and, as a consequence, a finite number of geometric factors. Using the equivalent relations between the Riemann mapping coefficients, the GPTs and the geometric factors,
one can determine $$C,\quad \{b_k\}_{k\leq N},\quad\{\sigma_k\}_{k\leq N-1},\quad\{\mu_k\}_{k\leq N-2}$$
from $$\left\{\gamma^2_{k1}\right\}_{k\leq N},\quad N\geq 2.$$
As a consequence, we have $\Phi_{N-2}[\Om]$ and $\Theta_{N-1}[\p\Om^r]$, where 
$\Phi_m[\Om]$ and $\Theta_m[\p \Om^r]$, $m\geq1$, are the truncation of $\Phi[\Om]$ and the Fourier series of $\Theta[\p \Om^r]$ at the $m$-th order, {\it i.e.},
\begin{align}
\Phi_m[\Om](\zeta)&=C\sum_{k=-1}^m \mu_k\zeta^{-k},\quad \zeta=\frac{1}{e^{2\pi t{\rm{i}}}},\label{eqn:PhiN}\\
\label{eqn:thetaN}
\Theta_m[\p\Om^r](t) &= 2 + 2\sum_{k=1}^m \Re\{\sigma_k\} \cos(2 \pi
kt) + 2\sum_{k=1}^m \Im\{\sigma_k\}\sin(2\pi kt)
\end{align}
for $t\in[0,1]$. 
From \eqnref{def:theta}, $\Theta_m[\p\Om^r](t)$ has an isolated peak at $t=t_0$ for a large number $m$ if $S[\Om^r](e^{2\pi t_0{\rm{i}}})$ is a corner point of $\p\Om^r$. For such $t_0$, $\Phi[\Om](e^{-2\pi t_0{\rm i}})$ is a corner point of $\p \Om$.

\section{Numerical results}\label{sec:experi}

\subsection{Description of the numerical method}
\label{sec:RCIP}

The numerical evaluation of the right-hand side of~(\ref{gpt})
involves, as its chief difficulty, the discretization and solution of
a Fredholm second kind integral equation on the piecewise smooth
boundary $\partial\Omega$. For this, we use {\it Nystr\"om
  discretization}~\cite[Chapter~4.1]{Atki97} based on 16-point
composite Gauss--Legendre quadrature and a computational mesh that is
dyadically refined in the direction of the corner vertices on
$\partial\Omega$. The resulting linear system for values of the
unknown layer density at the discretization points is compressed using
a lossless technique called {\it recursively compressed inverse
  preconditioning} (RCIP)~\cite{HelsOjal08} and solved using a
standard direct method. We have implemented our scheme in {\sc
  Matlab}. The execution time in the numerical examples in this paper
is typically a few seconds.

The RCIP technique serves two purposes. First, it greatly accelerates
the solution process when boundary singularities are present. In fact,
the combination of Nystr\"om discretization and RCIP acceleration
enables the solution of Fredholm second kind integral equations on
piecewise smooth boundaries with approximately the same speed at which they can be solved on smooth boundaries using Nystr\"om
discretization only. Second, RCIP stabilizes the solution process to
the extent that integral equations modeling well-conditioned boundary
value problems for elliptic partial differential equations in
piecewise smooth domains often can be solved with almost machine
precision. See~\cite{HKL17,HelsKarl16,HMM-NJP-11,HP-ACHA-13} for
examples where RCIP accelerated Nystr\"om discretization has been used
to compute, very accurately, polarizabilities and resonances of
various arrangements of dielectric objects with sharp corners and
edges. See also the recently revised compendium~\cite{Tutorial-arXiv}
for a comprehensive review of the RCIP technique and an ample
reference list.

It should be mentioned that RCIP compresses the integral equation
around one corner of $\partial\Omega$ at a time. The compression
requires, for each corner, a local boundary parameterization $z_{\rm
  loc}(t)$. If the original parameterization $z(t)$ has a corner
vertex at $t=t_i$, then this local parameterization is defined by
\begin{equation}
z_{\rm loc}(t):= z(t+t_i)-z(t_i)\,.
\label{eq:locparam}
\end{equation}
This means that mesh refinement occurs at $t\approx 0$ and that
$z_{\rm loc}(0)$ is at the origin. For high achievable accuracy in the
solution, the numerical implementation of $z_{\rm loc}(t)$ from the
definition~(\ref{eq:locparam}) is usually not good due to numerical
cancellation. Rather, $z_{\rm loc}(t)$ should be available in a form
that allows for evaluation with high relative accuracy also for small
arguments $t$. In the present work we find local parameterizations
accurate for small arguments using series expansion techniques.

\subsection{Symmetric domain}
\label{example:shape1}

\begin{figure}[p]
        \centering
    \begin{subfigure}[t]{0.45\textwidth}
        \centering
        \includegraphics[height = 4.5cm, width=4.5cm]{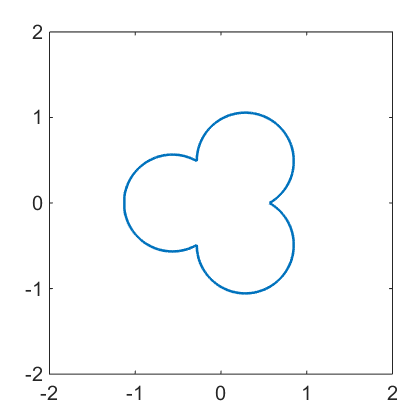}\caption{Symmetric curvilinear triangle $\Omega$}\label{symm_a}
    \end{subfigure}
\begin{subfigure}[t]{0.45\textwidth}
        \centering
        \includegraphics[height = 4.5cm, width=4.5cm]{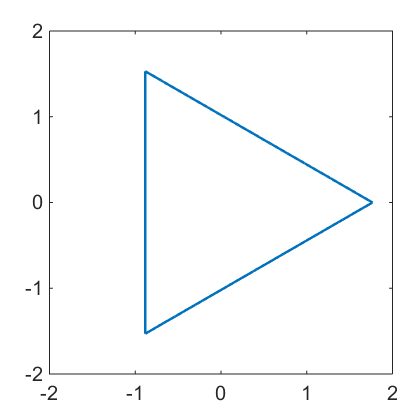}\caption{Equilateral triangle $P=\Om^r$}\label{symm_b}
    \end{subfigure}
    \\[2mm]
    \begin{subfigure}[t]{0.45\textwidth}
        \centering
        \includegraphics[height = 4.5cm, width=5.5cm]{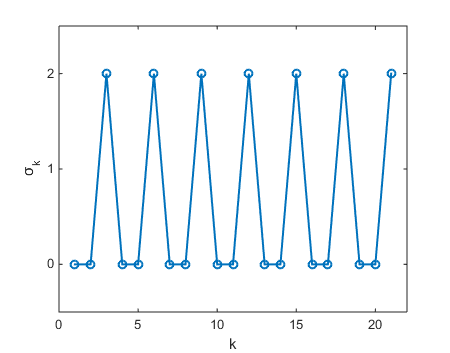}\caption{Geometric factor $\sigma_k\left[\Om\right]$}\label{symm_c}
    \end{subfigure}
    \begin{subfigure}[t]{0.45\textwidth}
        \centering
        \includegraphics[height = 4.5cm, width=6.5cm]{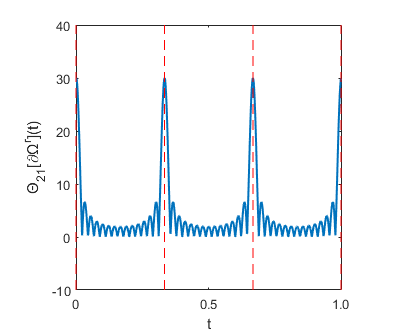}\caption{Truncation of $\Theta[\p\Om^r]$}\label{symm_d}
    \end{subfigure}
\caption{Symmetric cornered domain. (a) and (b) illustrate a
  curvilinear triangle $\Om$ and its reflection $P$ across the unit
  circle. Since ${\Om}$ is $3$-point radially symmetric,
  $\sigma_k[{\Om}] = 0$ for every $k \not \equiv 0 \pmod 3$ as shown in (c).
  {The graph $\Theta_{21}[\p\Om^r]$ in (d) shows three isolated peaks,
  located on the red dashed vertical lines indicating $t$-values for
  which $\Phi[\Om](e^{-2\pi t{\rm{i}}})$ are corner points of $\p\Om$.}}\label{figure:triangle}
\end{figure}

In this example, we consider a symmetric curvilinear triangle $\Om$.
The domain $\Om$ is the reflection of an equilateral triangle
$P$ across the unit circle. Note that $P=\Om^r$. See Figure
\ref{figure:triangle}(a,b) for the shape of $\Om$ and $P$.

The interior Riemann mapping function corresponding to the triangle
$P$ is
\begin{equation*}
S[\Om^r](z) = \frac{1}{C} \int_{0}^{z} \prod_{j=1}^3
\left(1 - \frac{w}{a_j}\right)^{-\frac{2}{3}} dw\,,
\end{equation*}
with $C=1$ and $a_j = \exp(\frac{2\pi j}{3} {\rm i})$ for $j=1,2,3$.
Since $P$ is a polygon, the boundary curvature is zero except at
corner points.  Hence, we have $\Theta[\p \Om^r](t) = \sum_{j=1}^3
\frac{2}{3}\delta(t-\frac{j}{3})$, and $\sigma_k[\Om]$ get periodic
values
\begin{equation}
\sigma_k\left[\Om\right]=
\begin{cases}
\ds 0,\quad&\mbox{if }k\not\equiv0 \quad (\mbox{mod }3)\,,\\[1mm]
\ds 2,\quad&\mbox{if }k\equiv0 \quad (\mbox{mod }3)\,.
\end{cases}\label{sigmak:example1}
\end{equation}
This fact fits well with Corollary \ref{lemma:sigmadecay} and the
existence of three corners on $\partial\Omega$. Since $\Om$ is
$3$-point radially symmetric, $\sigma_k[\Om] = 0$ for every $k \not
\equiv 0 \pmod 3$ as shown in Corollary \ref{radialsymm}.  

Figure
\ref{figure:triangle}(c) shows the graph of the geometric factors
and Figure \ref{figure:triangle}(d) the graph of $\Theta_{21}[\p\Om^r](t)$.
$\Theta_{21}[\p\Om^r](t)$ exhibits three isolated peaks at $t$-values corresponding to corner points, which are marked by red vertical dashed lines. The peaks correspond to the Dirac delta singularities in $\Theta$.

Now we perform a numerical computation to solve \eqnref{gpt} using the RCIP-accelerated
Nystr\"om scheme described in Section \ref{sec:RCIP} and acquire the
GPTs from \eqnref{GPT12}. Using the computed GPTs, we
then calculate $\sigma_k$ via~(\ref{eqn:bk1}) and~(\ref{sigma}). Table
\ref{table:sigma_ex1} displays the 20 first computed values of
$\sigma_k$. The acquired non-zero values agree with the analytic
values in \eqnref{sigmak:example1} to between 10 and 14 digits. The
zero values agree even better. 

\begin{table}[p]
{\begin{tabular}{|c|c|c|c|}
\hline
$k$ & $\sigma_k$ & $k$ & $\sigma_k$     \\
\hline
$1$ &  $0$                  & $11$ & $-1\times10^{-15}$\\
\hline
$2$ &  $-3\times10^{-17}$ & $12$ & $2.00000000001$\\
\hline
$3$ &  $1.999999999999994$  & $13$ & $1\times10^{-15}$\\
\hline
$4$ &  $1\times10^{-33}$    & $14$ & $1\times10^{-16}$\\
\hline
$5$ &  $-3\times10^{-16}$   & $15$ & $2.0000000001$\\
\hline
$6$ &  $1.99999999999997$   & $16$ & $8\times10^{-16}$\\
\hline
$7$ &  $2\times10^{-32}$    & $17$ & $-3\times10^{-15}$\\
\hline
$8$ &  $-3\times10^{-16}$   & $18$ & $2.000000001$\\
\hline
$9$ &  $2.0000000000003$    & $19$ & $-7\times10^{-16}$\\
\hline
$10$ & $4\times10^{-16}$    & $20$ & $-7\times10^{-15}$\\
\hline
\end{tabular}}
\caption{The 20 first geometric factors $\sigma_k$ of the symmetric
domain in Figure \ref{figure:triangle}(a), obtained
from numerically computed GPTs. The values agree well
with the analytic expression of \eqnref{sigmak:example1}.}
\label{table:sigma_ex1}
\end{table}

\subsection{Non-symmetric domain}
\label{example:shape2}

In this example $\Om$ is the non-symmetric domain with corners in
Figure \ref{fig:asymm}(a) (see Appendix \ref{appen:para} for the
parametrization). As in Section~\ref{example:shape1}, the GPTs are
numerically computed and the geometric factors $\sigma_k$ are
calculated from the GPTs via (\ref{GPT12}), (\ref{gpt}),
(\ref{eqn:bk1}) and~(\ref{sigma}). Table~\ref{table:sigma_ex2}
displays the first 20 geometric factors. {Note that
$\sigma_k\left[\Om\right]$ shows an oscillatory behavior as $k$ increases, as also shown in Figure
\ref{fig:asymm}(c).} 
The graph of $\Theta_{28}[\p \Om^r](t)$ shows three isolated
peaks at the locations of the corner points, which are marked by red
vertical lines. Again, the peaks correspond to the Dirac delta
singularities in $\Theta$.
\begin{figure}[p]
        \centering
    \begin{subfigure}[t]{0.4\textwidth}
        \centering
        \includegraphics[height = 4.5cm, width=4.5cm]{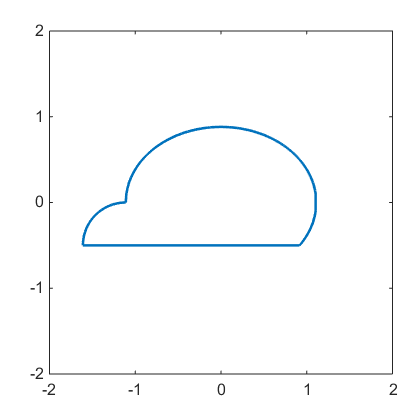}
        \caption{Cap-shaped cornered domain $\Om$}\label{asymm_a}
    \end{subfigure}
    \begin{subfigure}[t]{0.4\textwidth}
        \centering
        \includegraphics[height = 4.5cm, width=4.5cm]{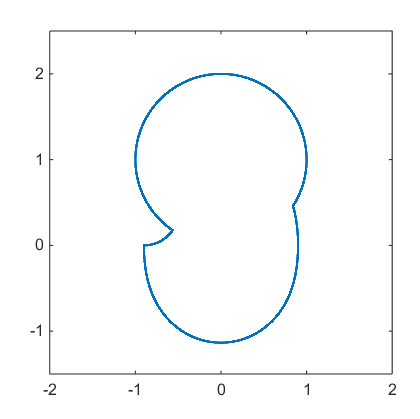}
        \caption{Reflected domain $\Om^r$}\label{asymm_b}
    \end{subfigure}
    \\
    \begin{subfigure}[t]{0.4\textwidth}
        \centering
        \includegraphics[height = 4.5cm, width=5.5cm]{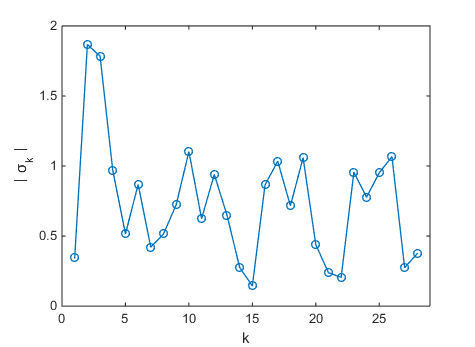}
        \caption{Geometric factor $\sigma_k[\Om]$}\label{asymm_c}
    \end{subfigure}
    \begin{subfigure}[t]{0.4\textwidth}
        \centering
        \includegraphics[height = 4.5cm,width=6.5cm]{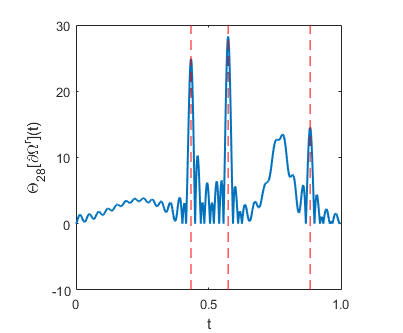}
        \caption{Truncation of $\Theta[\p \Om^r]$}\label{asymm_d}
    \end{subfigure}
\caption{Non-symmetric cornered domain. (a) and (b) illustrate a
  non-symmetric domain $\Om$ and its reflection $\Om^r$ across the
  unit circle. {(c) shows an oscillatory behavior of $\sigma_k[\Om]$ as $k$ increases, which is consistent with Corollary \ref{lemma:sigmadecay}.}  The graph $\Theta_{28}[\p\Om^r](t)$ in (d) shows three isolated
  peaks, located on the red dashed vertical lines indicating
  $t$-values for which $\Phi[\Om](e^{-2\pi t{\rm{i}}})$ are corner
  points of $\p\Om$}\label{fig:asymm}  
\end{figure}
\begin{table}[p]
{\begin{tabular}{|c|c|c|c|}
\hline
$k$ & $\sigma_k$ & $k$ & $\sigma_k$ \\
\hline
$1$  & $ 0.336144826114240 - 0.076400757440234{\rm i}$ & $11$ &
       $ 0.186721819078    - 0.595616065541{\rm i}$\\
\hline
$2$  & $-1.75172536453942  - 0.64675188584893{\rm i}$  & $12$ &
       $ 0.023499384397    + 0.939740685579{\rm i}$\\
\hline
$3$  & $ 0.03406793409600  + 1.78388113685821{\rm i}$  & $13$ &
       $ 0.49141111608     - 0.42493152167{\rm i}$\\
\hline
$4$  & $ 0.82403911365013  - 0.50742133234639{\rm i}$  & $14$ &
       $ 0.14107795747     - 0.23706422509{\rm i}$\\
\hline
$5$  & $ 0.49942961065065  + 0.12520990108117{\rm i}$  & $15$ &
       $-0.11870768404     + 0.08096200618{\rm i}$\\
\hline
$6$  & $-0.1083287142652   - 0.8609605708526{\rm i}$   & $16$ &
       $ 0.0444798323      - 0.8639904912{\rm i}$\\
\hline
$7$  & $-0.3918884064658   - 0.1538531930618{\rm i}$  & $17$ &
       $-0.6290121604      + 0.8209980128{\rm i}$\\
\hline
$8$  & $-0.147595479311    + 0.493831447944{\rm i}$   & $18$ &
       $ 0.1981029905      - 0.6878495747{\rm i}$\\
\hline
$9$  & $-0.072598481664    - 0.719868325959{\rm i}$   & $19$ &
       $-0.4512400133      + 0.9581338937{\rm i}$\\
\hline
$10$ & $-0.330200317533    + 1.052309941949{\rm i}$   & $20$ &
       $ 0.407788339       - 0.162013757{\rm i}$\\
\hline
\end{tabular}}
\caption{The 20 first geometric factors $\sigma_k$ of the
  non-symmetric domain in Figure
  \ref{fig:asymm}(a), obtained from numerically computed GPTs.}
\label{table:sigma_ex2}
\end{table}

\begin{figure}[h]
        \centering
             \includegraphics[height = 4.5cm,width=5.5cm]{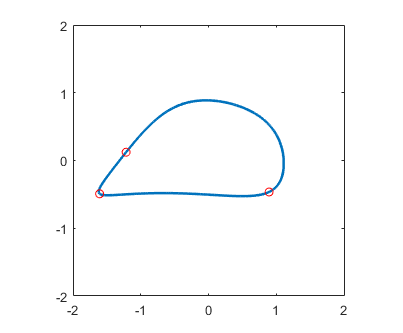}\hskip .6cm
              \includegraphics[height = 4.5cm, width=6.5cm]{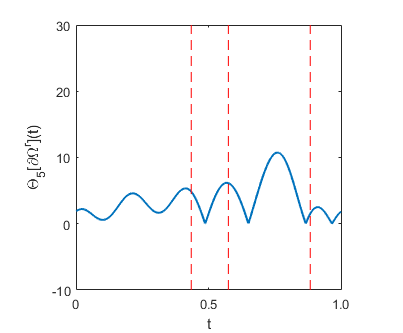}\\[.5cm]
              \includegraphics[height = 4.5cm,width=5.5cm]{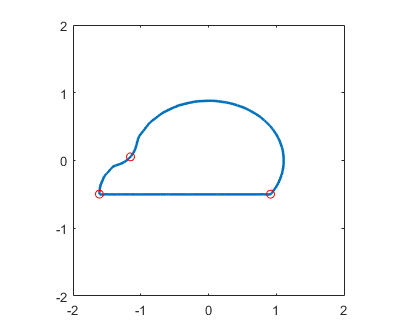}\hskip .6cm
             \includegraphics[height = 4.5cm, width=6.5cm]{Fourier_peak2}
           \caption{Imaging $\Om$ from a finite number of components of the GPTs. $\Om$ is given as in Figure \ref{fig:asymm}(a) and $\left\{\gamma^2_{k1}: k\leq N\right\}$ are used with $N=6$ (top) and $N=29$ (bottom). The left column shows the images of the unit circle under $\Phi_{N-2}[\Om]$. {The small red circles indicate values of $\Phi_{N-2}[\Om](e^{-2\pi
  t_0{\rm{i}}})$ for which $\Phi[\Om](e^{-2\pi t_0 {\rm{i}}})$ are
  corner points of $\p\Om$.} The right column illustrates $\Theta_{N-1}[\p\Om^r]$. For small $N$ there is no isolated peak for $\Theta_{N-1}[\p\Om^r]$. }
\label{fig:truncation}
\end{figure}

In Figure \ref{fig:truncation}, we consider the imaging problem of $\Om$ from a finite number of components of the GPTs. As discussed in Section \ref{sec:finiteGPTs} one can have the truncated series $\Phi_{N-2}[\Om](\zeta)$ from  $\left\{\gamma^2_{k1}: k\leq N\right\}$. In this example, we reconstruct $\Om$ using $N=6$ and $N=29$.
The image of the unit circle under $\Phi_{N-2}[\Om]$ approximates the shape of $\p \Om$ even for small $N$. For $N=29$ the graph of $\Theta_{N-1}[\p\Om^r](t)$ shows isolated peaks at $t$ corresponding to corner points, while  it does not for $N=6$.
%
%
%

\subsection{Smooth domain}
\begin{figure}[h]
        \centering
     \begin{subfigure}[t]{0.4\textwidth}
        \centering
        \includegraphics[height = 4.5cm, width=4.5cm]{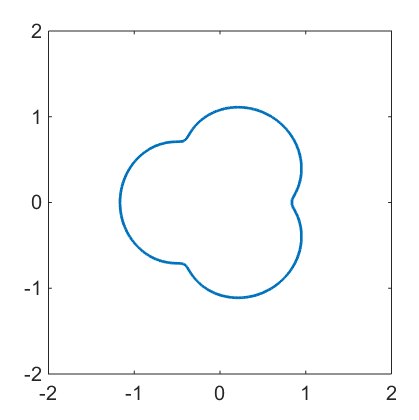}\caption{Symmetric smooth domain $\widetilde{\Omega}$}\label{symm_smooth_a}
    \end{subfigure}
    \begin{subfigure}[t]{0.4\textwidth}
        \centering
        \includegraphics[height = 4.5cm, width=4.5cm]{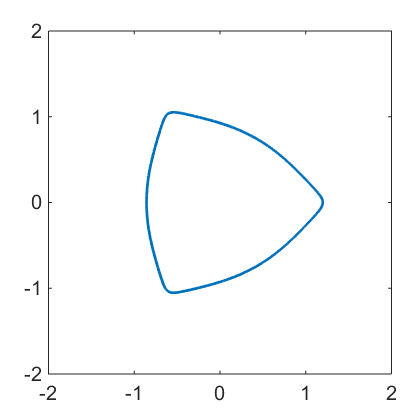}\caption{Reflected domain ${\widetilde{\Om}}^r$}\label{symm_smooth_b}
    \end{subfigure}
     \\
     \begin{subfigure}[t]{0.4\textwidth}
        \centering
        \includegraphics[height = 4.5cm, width=5.5cm]{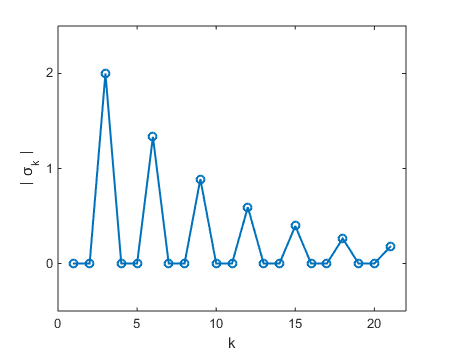}\caption{Geometric factor $\sigma_k[\widetilde{\Om}]$}\label{symm_smooth_c}
    \end{subfigure}
    \begin{subfigure}[t]{0.4\textwidth}
        \centering
        \includegraphics[height = 4.5cm, width=6.5cm]{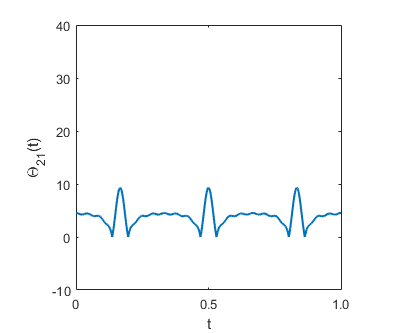}\caption{Truncation of $\Theta[\p\widetilde{\Om}^r]$}\label{symm_smooth_d}
    \end{subfigure}
\caption{Symmetric smooth domain. (a) and (b) illustrate a smooth
  domain $\widetilde{\Om}$ and its reflection $\widetilde{\Om}^r$
  across the unit circle. (c) shows the exponential decay of $\sigma_k$ as $k$ increases. Since $\widetilde{\Om}$ is $3$-point
  radially symmetric, $\sigma_k[\widetilde{\Om}] = 0$ for every $k
  \not \equiv 0 \pmod 3$.
 (d) shows the graph of a truncated Fourier series of $\Theta[\p \Om^r]$. $\Theta_{21}[\p
  \widetilde{\Om}^r]$ does not exhibit isolated peaks, but has
  relatively large values at boundary points of $\widetilde{\Om}^r$
  with large curvature.}\label{figure:smoothtriangle}
\end{figure}

\begin{figure}[h]
        \centering
     \begin{subfigure}[t]{0.4\textwidth}
        \centering
        \includegraphics[height = 4.5cm, width=4.5cm]{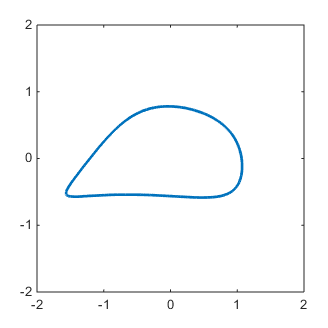}
        \caption{Symmetric smooth domain $\widetilde{\Omega}$}\label{asymm_smooth_a}
    \end{subfigure}
    \begin{subfigure}[t]{0.4\textwidth}
        \centering
        \includegraphics[height = 4.5cm, width=4.5cm]{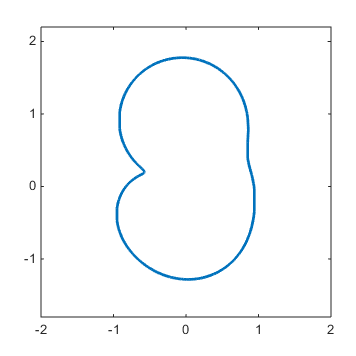}
        \caption{Reflected domain ${\widetilde{\Om}}^r$}\label{asymm_smooth_b}
    \end{subfigure}
     \\
     \begin{subfigure}[t]{0.4\textwidth}
        \centering
        \includegraphics[height = 4.5cm, width=5.5cm]{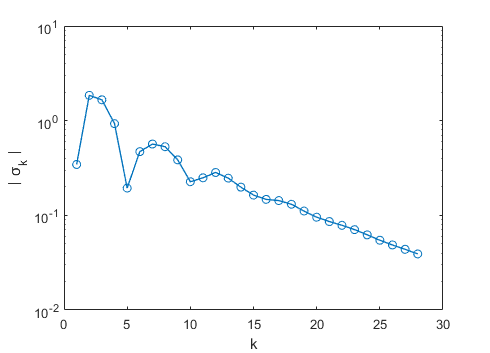}
        \caption{Graph of $\sigma_k[\widetilde{\Om}]$ in log scale}\label{asymm_smooth_c}
    \end{subfigure}
    \begin{subfigure}[t]{0.4\textwidth}
        \centering
        \includegraphics[height = 4.5cm, width=6.5cm]{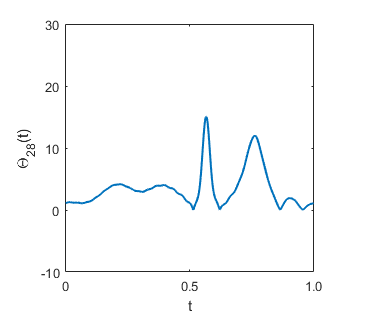}
        \caption{Truncation of $\Theta[\p\widetilde{\Om}^r]$}\label{asymm_smooth_d}
    \end{subfigure}
\caption{Non-symmetric smooth domain. (a) and (b) illustrate a smooth
  domain $\widetilde{\Om}$ and its reflection $\widetilde{\Om}^r$
  across the unit circle. (c) shows the exponential decay of $\sigma_k$ as $k$ increases. 
  The graph of $\Theta_{28}[\p  \widetilde{\Om}^r]$ in (d) does not exhibit isolated peaks, but has
  relatively large values at boundary points of $\widetilde{\Om}^r$
  with large curvature.}
\label{fig:asymmetricsmooth}
\end{figure}

In Figure \ref{figure:smoothtriangle} we consider a smooth domain, denoted by
$\widetilde{\Om}$, with $3$-point radial symmetry. Note that
$\widetilde{\Om}$ has a similar shape as that of $\Om$ in Figure
\ref{figure:triangle}. In Figure \ref{fig:asymmetricsmooth} we
consider another smooth domain, denoted again by $\widetilde{\Om}$, which has a similar
shape as $\Om$ in Figure \ref{fig:asymm}.

In both cases the domain $\widetilde{\Om}$ is made as
$\p\widetilde{\Om}=\{P(z):|z|=1\}$ by using a polynomial $P(z)$. Then,
the corresponding geometric factors are analytically calculated. In
contrast to the examples with the cornered domain, the geometric factors
of $\widetilde{\Om}$ decay exponentially. The partial sums of Fourier series of $\Theta$
have relatively large values at $t$-values corresponding to boundary points of $\widetilde{\Om}^r$
with large curvature.

Since $\widetilde{\Om}$ in Figure \ref{figure:smoothtriangle} is
$3$-point radially symmetric, we have $\sigma_k[\widetilde{\Om}] = 0$
for every $k \not \equiv 0 \pmod 3$.

\section{Proofs}
\label{sec:proof}

\subsection{Proof of Proposition \ref{theorem1}}\label{proof:theorem1}

By \eqnref{Phi}, \eqnref{RG:conformal}, and \eqnref{Sseries0} and the
definitions of $\mu_{n,k}$ and $b_{n,k}$, we have 
\begin{align} \quad
 & \Phi(\zeta)^n - \sum_{m=1}^\infty \frac{ \gamma_{mn}^1
    }{\Phi(\zeta)^m}
  = \Phi(\zeta)^n - \sum_{m=1}^\infty \gamma_{mn}^1  S\Bigr(\frac{1}{\zeta}\Bigr)^m \nonumber \\
  &= C^n\zeta^{2n}\Bigr(\sum_{k=1}^\infty \frac{\mu_{k-2}}{\zeta^k} \Bigr)^n- \sum_{m=1}^{\infty}  \gamma_{mn}^1 \frac{1}{C^m} \Bigr( \sum_{k=1}^\infty \frac{b_{k}}{\zeta^k} \Bigr)^m \nonumber
  = C^n \sum_{k=n}^{\infty} \frac{\mu_{n,k}}{\zeta^{k-2n}}- \sum_{k=1}^{\infty}\Bigr(\sum_{m=1}^{k}\frac{1}{C^m} \gamma_{mn}^1 b_{m,k} \Bigr)\frac{1}{\zeta^k} \nonumber \\
  &= C^n \sum_{k=0}^{n} \mu_{n,2n-k}\zeta^{k}+\sum_{k=1}^{\infty}
  \Bigr( C^n \mu_{n,2n+k}-\sum_{m=1}^{k} \frac{1}{C^m} \gamma_{mn}^1
  b_{m,k} \Bigr)\frac{1}{\zeta^k}.  \label{rel}
\end{align}
{Note that $$ \Phi(\zeta)^n - \sum_{m=1}^\infty \frac{ \gamma_{mn}^1    }{\Phi(\zeta)^m}=\frac{1}{2}\left(V_1\circ\Phi+V_2\circ\Phi\right)(\zeta)$$
with $V_1$ and $V_2$ defined in Section \ref{sec:GPTandcoeff}. 
From Lemma \ref{lemma:basic}, the right-hand side of the above equation has an entire extension. The main idea of the proof of Lemma \ref{lemma:basic} is the equality
\begin{align}
 -\frac{1}{2}{\left(V_1\circ\Phi+V_2\circ\Phi\right)(\zeta)}\label{eqn:V1V2}
&=\mbox{const.}-\frac{1}{2}\overline{\left(V_1\circ\Phi-V_2\circ\Phi\right)(\zeta)}\\
&=\mbox{const.}-\frac{1}{2}\overline{\left(V_1\circ\Phi-V_2\circ\Phi\right)(\frac{1}{\bar{\zeta}})}\notag
\end{align}
on $|\zeta|=1$ due to \eqnref{extension1} and \eqnref{extension2}. 
Since the last term in the above equation is analytic for $|\zeta|<1$, the function $\frac{1}{2}{\left(V_1\circ\Phi+V_2\circ\Phi\right)(\zeta)}$ has an entire extension. 
Therefore, the principal parts of \eqnref{rel}
should vanish: \beq\label{gamma1_all} C^n
\mu_{n,2n+k}-\sum_{m=1}^{k} \frac{1}{C^m} \gamma_{mn}^1 b_{m,k}=0.
\eeq Rearranging \eqnref{gamma1_all}, and since $b_{k,k}= 1$, we get 
$
\gamma_{kn}^1 = C^{k+n} \Bigr(\mu_{n,2n+k}-\sum_{m=1}^{k-1} \frac{1}{C^{m+n}} \gamma_{mn}^1 b_{m,k} \Bigr)$ for each $ n,k\in \mathbb{N}$.
%
This proves \eqnref{Gamma1}.

From \eqnref{eqn:V1V2}, for each $k\in\mathbb{N}\cup\{0\}$ and sufficiently large $R>0$ we have
\begin{align}
&\int_{|\zeta|=1}-\frac{1}{2}{\left(V_1\circ\Phi-V_2\circ\Phi\right)(\zeta)}\,\zeta^k d\zeta\notag
= \int_{|\zeta|=1}\overline{-\frac{1}{2}{\left(V_1\circ\Phi+V_2\circ\Phi\right)(\zeta)}}\,\zeta^k d\zeta\\\notag
&=\int_{|\zeta|=1}\overline{-\frac{1}{2}{\left(V_1\circ\Phi+V_2\circ\Phi\right)(\bar{\zeta}^{-1})}}\,\zeta^k d\zeta=\int_{|\zeta|=\frac{1}{R}}\overline{-\frac{1}{2}{\left(V_1\circ\Phi+V_2\circ\Phi\right)(\bar{\zeta}^{-1})}}\,\zeta^k d\zeta\\\label{eqn:formal1}
&=-\int_{|\zeta|=\frac{1}{R}}\overline{\Phi(\bar{\zeta}^{-1})^n - \sum_{m=1}^\infty \frac{ \gamma_{mn}^1    }{\Phi(\bar{\zeta}^{-1})^m}}\,\zeta^k d\zeta,
\end{align}
and the left-hand side equals
\beq\label{eqn:formal2}
\int_{|\zeta|=1}-\frac{1}{2}{\left(V_1\circ\Phi-V_2\circ\Phi\right)(\zeta)}\,\zeta^k d\zeta=
\int_{|\zeta|=R}\left(\sum_{m=1}^\infty \frac{\gamma_{mn}^2}{\Phi(\zeta)^m}\right)\zeta^kd\zeta.
\eeq

We compute
\begin{align}
\overline{\Phi(\bar{\zeta}^{-1})^n - \sum_{m=1}^\infty \frac{ \gamma_{mn}^1    }{\Phi(\bar{\zeta}^{-1})^m}}
&=\overline{\Phi(\bar{\zeta}^{-1})^n} - \overline{\sum_{m=1}^\infty { \gamma_{mn}^1    }{S(\bar{\zeta})^m}}\notag\\
&={C}^n\Bigr(\sum_{k=-1}^\infty \overline{\mu_k}\zeta^k \Bigr)^n-\sum_{m=1}^{\infty}  \overline{\gamma_{mn}^1}\frac{1}{{C}^m} \Bigr(\sum_{k=1}^\infty \overline{b_k}\zeta^k \Bigr)^m \quad \mbox{on }|\zeta|=\frac{1}{R},\label{eqn:formal3}
\end{align}
and
\beq\label{eqn:formal4}\sum_{m=1}^\infty \frac{\gamma_{mn}^2}{\Phi(\zeta)^m}=\sum_{m=1}^\infty \gamma_{mn}^2S(\frac{1}{\zeta})^m=\sum_{m=1}^{\infty}  \gamma_{mn}^2 \frac{1}{C^m} \Bigr(\sum_{k=1}^\infty \frac{b_{k}}{\zeta^k}\Bigr)^m\quad \mbox{on }|\zeta|=R.\eeq
By applying a similar multinomial expansion as in \eqnref{rel}, we can formally expand the summation of two components in \eqnref{eqn:formal3} and \eqnref{eqn:formal4} which contain principal parts:}
\begin{align}
&\sum_{m=1}^{\infty}  \gamma_{mn}^2\frac{1}{C^m} \Bigr(\sum_{k=1}^\infty \frac{b_{k}}{\zeta^k}\Bigr)^m+{C}^n\Bigr(\sum_{k=-1}^\infty \overline{\mu_k}\zeta^k \Bigr)^n \nonumber \\
= &\sum_{k=1}^{\infty} \Bigr(\sum_{m=1}^{k} \frac{1}{C^m} \gamma_{mn}^2 b_{m,k} \Bigr)\frac{1}{\zeta^k}+ {C}^n \sum_{k=n}^{\infty}\overline{\mu_{n,k}}\zeta^{k-2n} \nonumber \\
=& {C}^n \sum_{k=0}^{\infty}\overline{\mu_{n,2n+k}}\zeta^{k} + \sum_{k=1}^{n} \Bigr({C}^n \overline{\mu_{n,2n-k}}+\sum_{m=1}^{k} \frac{1}{C^m} \gamma_{mn}^2 b_{m,k}\Bigr)\frac{1}{\zeta^{k}} + \sum_{k=n+1}^{\infty} \Bigr(\sum_{m=1}^{k} \frac{1}{C^m} \gamma_{mn}^2 b_{m,k}\Bigr)\frac{1}{\zeta^{k}}.        \notag
\end{align}
In view of \eqnref{eqn:formal1} and \eqnref{eqn:formal2}, the principal
parts of the above equation should vanish:
\beq\notag
\begin{cases}
 {C}^n \overline{\mu_{n,2n-k}} + \sum_{m=1}^{k} \frac{1}{C^m} \gamma_{mn}^2 b_{m,k} = 0, &\mbox{for }1\le k\le n,     \\
 \sum_{m=1}^{k} \frac{1}{C^m} \gamma_{mn}^2 b_{m,k} = 0, &\mbox{for }k\ge n+1.         
 \end{cases}
\eeq
Rearranging the above equation, and since
$b_{k,k}= 1$, we get
\beq\notag
\gamma_{kn}^2
= 
\begin{cases}
- C^{k+n}\Bigr( \overline{\mu_{n,2n-k}} + \sum_{m=1}^{k-1}\frac{1}{C^{m +n}}\gamma_{mn}^2 b_{m,k}\Bigr), &\mbox{for }1\le k\le n,     \\
- C^{k+n} \sum_{m=1}^{k-1} \frac{1}{C^{m+n}} \gamma_{mn}^2 b_{m,k}, &\mbox{for }k\ge n+1,
\end{cases}
\eeq
This proves \eqnref{Gamma2}.
\qed

\subsection{Proofs of  Proposition \ref{GPTtob} and  Remark \ref{remark:bk2}}
\label{sec:proof:others}

\noindent{\textbf{Proof of Proposition \ref{GPTtob}}}
From \eqnref{Gamma2} with $n=k=1$ and \eqnref{bkk}, we have
$\gamma_{11}^2 = -C^2 \overline{\mu_{1,1}} = -C^2.$
This implies \eqnref{eqn:C}.
Applying again \eqnref{Gamma2} for $n=1$ and $k\ge2$, we have
$
\gamma_{k1}^2 = -C^{k+1} \sum_{m=1}^{k-1} \frac{\gamma_{m1}^2}{C^{m+1}}  b_{m,k},
$
which is equivalent to $\sum_{m=1}^{k} \frac{\gamma_{m1}^2}{C^{m+1}}  b_{m,k} = 0$.
Hence we have for $k\geq2$,
\beq
0 = \sum_{m=1}^{k} \frac{\gamma_{m1}^2}{C^{m+1}}  b_{m,k} = \frac{\gamma_{11}^2}{C^{2}}  b_{1,k} + \sum_{m=2}^{k} \frac{\gamma_{m1}^2}{C^{m+1}}  b_{m,k} = -b_k + \sum_{m=2}^{k} \frac{\gamma_{m1}^2}{C^{m+1}}  b_{m,k}\,.\nonumber
\eeq
This proves \eqnref{eqn:bk1}.
\qed

\vskip .3cm
\noindent{\textbf{Proof of Remark \ref{remark:bk2}}}
From \eqnref{Gamma1}, we have
\begin{align*}
\gamma_{k1}^1= C^{k+1} \Bigr(\mu_{1,2+k}-\sum_{m=1}^{k-1} \frac{\gamma_{m1}^1}{C^{m+1}}  b_{m,k} \Bigr), \quad \mbox{for } k\ge 2.
\end{align*}
It follows that $\mu_{1,2+k}=\sum_{m=1}^{k-1}\frac{\gamma_{m1}^1}{C^{m+1}}b_{m,k}+\frac{\gamma_{k1}^1}{C^{k+1}}=\sum_{m=1}^{k} \frac{\gamma_{m1}^1}{C^{m+1}}  b_{m,k}$ because $b_{k,k}=1$.
Hence, for $k\ge 2$, we have
\beq
\mu_{k} = \mu_{1,2+k}=\sum_{m=1}^{k} \frac{\gamma_{m1}^1}{C^{m+1}}  b_{m,k}= \frac{\gamma_{11}^1}{C^{2}} b_k + \sum_{m=2}^{k} \frac{\gamma_{m1}^1}{C^{m+1}}  b_{m,k}\,,\nonumber
\eeq
where $C$ is as in \eqnref{eqn:C}. This proves \eqnref{eqn:bk2}.
\qed

\section{Conclusions}

We have analyzed the effects of corners of an insulating inclusion
$\Om$ on the perturbation of an electric potential. We
derived explicit connections between generalized polarization tensors
and coefficients of interior Riemann mapping functions on the way. We defined a
sequence of geometric factors using these mapping coefficients. Mutually equivalent relations are then deduced  between GPTs, the Riemann
mapping coefficients, and the geometric factors. 
We finally characterized the corner effect: the sequence of geometric factors is the sequence of Fourier coefficients of the generalized external angle function for $\Om^r$, the reflection of $\Om$ across the unit circle, where the generalized external angle function contains the Dirac delta singularity at the corner points. Based on this corner effect, we established a criteria for the existence of corner points on the inclusion boundary in terms of the geometric factors. We assumed that the
inclusion is insulated. It will be of interest to find geometric
factors for inclusions with arbitrary conductivity that reveal the
presence of corners.

\begin{appendices}

\section{Carath\'eodory's mapping theorem}\label{sec:cara}
For the relation between the convergence of domains and the convergence of the corresponding conformal mappings, we introduce some content from \cite{markushevich1977theory}.

\begin{definition}
Let $\{ \Omega_n \}_{n\in \Nbb}$ be a sequence of simply connected and uniformly bounded domains in $\CC$, and each $\Omega_n$ contains a fixed disk centered at $z_0$. The ${kernel}$ of $\{ \Omega_n \}_{n\in \Nbb}$ is defined as the largest open domain $\Omega_{z_0}$ containing $z_0$ such that every compact subset $K \subset \Omega_{z_0}$ belongs to $\Omega_n$ for all $n \ge N$ with some $N \in \Nbb$ depending on $K$.
\end{definition}
\begin{definition}[Kernel convergence in the sense of Carath\'eodory]
Let $G_{z_0}$ be a kernel of $\{ \Omega_n \}_{n\in \Nbb}$ relative to the point $z_0$. If every subsequence of $\{ \Omega_n \}_{n\in \Nbb}$ has the same kernel $\Omega_{z_0}$, then $\Omega_n$ is said to converge to $\Omega_{z_0}$. Otherwise, $\Omega_n$ is said to diverge.
%
\end{definition}

We defined a concept of convergence in the sense of Carath\'eodory. Now, let's see how the convergence in the sense of Carath\'eodory is related to the convergence of the function sequences.


\begin{theorem}[Carath\'eodory's mapping theorem]
For each $n\in \Nbb$, let $f_n:\Omega_n \to \Dbb$ be a conformal mapping that satisfies
\beq
f_n(z_0) = 0, \quad f'_n(z_0)>0.
\nonumber
\eeq
Similarly, let $f:\Omega_{z_0} \to \Dbb$ be a conformal mapping that satisfies
\beq
f(z_0) = 0, \quad f'(z_0)>0.
\nonumber
\eeq
If $\Omega_n$ converges to $\Omega_{z_0}$, then $f_n$ converges uniformly to $f$ inside $\Omega_{z_0}$ (which means by definition that $f_n$ converges uniformly on any compact subset of $\Omega_{z_0}$), and $f_n^{-1}$ converges uniformly to $f^{-1}$ inside $\Dbb$.
Conversely, if $f_n$ converges uniformly to $f$ inside $\Omega_{z_0}$, or if $f_n^{-1}$ converges uniformly to $f^{-1}$ inside $\Dbb$, then $\Omega_n$ converges to $\Omega_{z_0}$.
\end{theorem}

\section{Parametrization of the non-symmetric domain in Section \ref{example:shape2}}\label{appen:para}
The boundary of the non-symmetric domain $\Om$ in Section \ref{example:shape2} can be parametrized as follows: 
\beq\notag
\gamma(t)=
\begin{cases}
\left(-\frac{1}{2}\sin\left(4\pi c t \right) - \frac{\sqrt{2}\pi}{4} , \quad - \frac{1}{2} + \frac{1}{2}\cos\left(4\pi c t \right)\right),\quad t\in[0,t_1],
\\[2mm]
\left(2\pi c(t-t_2) - \sqrt{2}\arcsin\left( \frac{\sqrt{2}}{2}\cos(2\pi a) \right) , \quad -\frac{1}{2} \right),\quad t\in[t_1,t_2],
\\[2mm]
\left(-\sqrt{2}\arcsin\left( \frac{\sqrt{2}}{2}\cos\left(2\pi c(t-t_2)+2\pi a\right) \right) , \quad -\arcsinh \left(\sin \left(2\pi c(t-t_2)+2\pi a\right)\right) \right), \quad t\in[t_2,1],
\end{cases}\eeq
with
$a = \frac{1}{2} - \frac{1}{2\pi}\arcsin\left(\sinh\left(\frac{1}{2}\right)\right)$, $b = a-\frac{1}{4\pi}-\frac{\sqrt{2}}{8}+\frac{\sqrt{2}}{2\pi}\arcsin\left(\frac{\sqrt{2}}{2}\cos(2\pi a)\right)$,
$c = \frac{9}{8}-b$, $t_1 = \frac{1}{8c} \approx 0.1122$, and $t_2 = t_1 + \left( \frac{a-b}{c} \right) \approx 0.4731.$

\end{appendices}

\bibliographystyle{plain}

{

}
\end{document}